\documentstyle[12pt]{article}

\input BoxedEPS.tex
\SetRokickiEPSFSpecial
\HideDisplacementBoxes
\newtheorem{thm}{Theorem}[section]
\newtheorem{lem}[thm]{Lemma}
\newtheorem{cor}[thm]{Corollary}
\newtheorem{rem}[thm]{Remark}
\newtheorem{defn}[thm]{Definition}
\newtheorem{prop}[thm]{Proposition}
\newtheorem{ex}[thm]{Example}

\begin{document}

\title{An algebraic annulus theorem}
\author{G. P. Scott and G. A. Swarup}
\maketitle

\begin{abstract}
We present an extension of Dunwoody's theory of tracks and use it to prove
an analogue of the annulus theorem for hyperbolic groups.
\end{abstract}

\tableofcontents

\section{ Introduction}

We prove the following analogue for (word) hyperbolic groups of the Annulus
Theorem for 3-manifolds:

\begin{thm}
Let $G$ be a torsion free hyperbolic group with one end. Suppose that $G$
has an infinite cyclic subgroup such that the number of ends of the pair $%
(G,H)$ is greater than one. Then $G$ splits over some infinite cyclic
subgroup.
\end{thm}

\begin{cor}
Let $G^{\prime }$ be a subgroup of finite index in a torsion free hyperbolic
group $G$. Then $G$ splits over an infinite cyclic group if and only if $%
G^{\prime }$ splits over an infinite cyclic group.
\end{cor}

This result has also been proved by Bowditch~\cite{bo:cut}, using very
different methods. We discuss these differences at the end of this
introduction.

The terminology used is standard (from~\cite{gh:book} and ~\cite{sc3:ends}).
The importance of splitting groups along infinite cyclic subgroups is well
known from the work of Paulin~\cite{pau:outer}, Rips and Sela~\cite{rs:rigid}
and Sela~\cite{sela:rigid2}.

If $G$ has an infinite cyclic subgroup $H$ such that {\em $e(G,H)\geq 3$, }%
we will say that $G$ is of {\it multi-band} {\it type}. Otherwise we will
say that $G$ is of {\it surface\ type.} The proofs of Theorem 1.1 are
different in the two cases. We do not need the assumption that $G$ is
torsion free in the multi-band case. In this case we obtain the following
result.

\begin{thm}
Let $G$ be a one-ended hyperbolic group and let $H$ be a two-ended subgroup
with $e(G,H)\geq 3$. Then $G$ splits over a subgroup commensurable with $H$.
\end{thm}

Both of the above results are closely related to the Annulus Theorem for
3-manifolds. There are several versions of this result and we state one of
the most basic ones here. A map $f:S^1\times I\rightarrow M$ of the annulus
into a 3-manifold $M$ is {\it essential} if it is $\pi _1$-injective, proper
in the sense that $f$ maps the boundary of the annulus into the boundary $%
\partial M$ of $M,$ and in addition $f$ is not properly homotopic into $%
\partial M.$

{\bf Annulus Theorem:}{\em \ Let }$M${\em \ be a compact orientable
irreducible 3-manifold with incompressible boundary. If }$M${\em \ admits an
essential map of the annulus }$S^1\times I,${\em \ then it admits an
essential embedding of the annulus.}

The connections between the Annulus Theorem for 3-manifolds and our
algebraic analogues are described briefly as follows. The assumptions that $%
M $ be compact orientable and irreducible with non-empty boundary imply that 
$G=\pi _1(M)$ is torsion free. The additional assumption that $M$ has
incompressible boundary implies that $G$ has one end. Also if $f:S^1\times
I\rightarrow M$ is a $\pi _1$-injective proper map and $H$ denotes the
infinite cyclic subgroup of $G$ carried by $f,$ then $f$ is essential if and
only if $e(G,H)\geq 2,$ and $G$ splits over $H$ if and only if $f$ is
properly homotopic to an embedding. Finally, if $e(G,H)\geq 3,$ then $f$ can
be homotoped to cover an embedded essential annulus.

Our paper is organised as follows. In section 2, we develop the general
theory of patterns and tracks in 2-complexes and its connection with the
number of ends of a 2-complex. This theory was introduced by Dunwoody in 
\cite{mjd:acc2}. We extend Dunwoody's theory by adding the idea of a
singular pattern and of the length of a pattern. Dunwoody introduced
patterns and tracks in order to prove the accessibility of finitely
presented groups. We show how our new ideas give an alternative version of
his arguments. In section 3, we discuss the use of tracks in a situation
where one has a pair of groups such that $e(G,H)\geq 2,$ and $H$ is infinite
cyclic. In section 4, we discuss crossing of tracks and give examples of
groups of multi-band type and of surface type. The proof of Theorem 1.3 on
the multi-band case is completed in section 5. The rest of the paper is
devoted to the surface type case, and the main construction of the proof
comes in section 8. We use hyperbolicity of the group a bit more in the
surface type case. On the whole our use of hyperbolicity seems weak and we
expect that the 2-complex techniques set out in our proofs will be useful in
the study of splittings of finitely presented groups. Another feature of the
proof is the surface-like properties of a large class of one-ended torsion
free hyperbolic groups, particularly those which we call surface type
groups. In this case, we need to deal with orientation reversing elements in
the group; we show that if a surface-like group has orientation reversing
elements then it splits as an amalgamated free product over an infinite
cyclic subgroup which has index two in one of the vertex groups which is
also infinite cyclic. This is the analogue of the fact that a non-orientable
surface must contain an embedded Mobius band. Our proof is suggested by the
arguments of Dunwoody~\cite{mjd:acc2} and Tukia~\cite{tukia:conj} and by the
least area arguments of Freedman, Hass and Scott~\cite{fhs:least} and \cite
{fhs:shortest}. As in Tukia~\cite{tukia:conj}, we run into difficulties in
the torsion case. With some more work one can show without the torsion free
assumption that either the conclusion of Theorem 1.1 holds or there is a
subgroup of $G$ which ``looks like'' a triangle group.

Brian Bowditch~\cite{bo:cut} recently developed a theory of
JSJ-decompositions for one-ended hyperbolic groups with locally connected
boundary, and Swarup then showed that any one-ended hyperbolic group has
locally connected boundary\cite{sw:cutpoint}. Using an extension of the work
of Rips, Bestvina and Feighn (see ~\cite{bo:tree}), Bowditch showed that if
one removes the torsion freeness hypothesis from Theorem 1.1, the conclusion
remains true except for the case of triangle groups. His techniques are very
different from ours, making heavy use of the boundary of a hyperbolic group,
and as explained above we expect that our techniques will be useful in the
study of splittings of finitely presented groups. It is for this reason that
we have presented the theory of patterns in more generality and detail than
we require in this paper. Perhaps Gabai's techniques~\cite{gabai:co} may
also be useful for extending Theorem 1.1. Another obvious question is
whether there is an analogue of orientation covers for hyperbolic groups.

Acknowledgments: Much of this work was carried out while the first author
was visiting the University of Melbourne in 1994 and the Mathematical
Sciences Research Institute in Berkeley in 1994 and 1995. He is grateful for
the partial support provided by the University of Melbourne, by MSRI and by
NSF grant DMS-9306240.

\section{Tracks and ends}

Throughout this section, $Y$ will denote a connected locally finite
2-dimensional simplicial complex. A subset of $Y$ will be called a {\it %
pattern} if it intersects each closed simplex of $Y$ in a compact properly
embedded submanifold of codimension one. This means that it avoids the
0-skeleton of $Y,$ meets each 1-simplex in a finite set and meets each
2-simplex in a compact 1-dimensional submanifold whose boundary lies on the
boundary of the simplex. See Figure 1a. A pattern $t$ in $Y$ locally has a
collar neighbourhood, and globally there is a neighbourhood of $t$ in $Y$
which is an $I$-bundle over $t$. If this bundle is trivial, we say that $t$
is {\it two-sided}. In this paper, we will mostly consider two-sided
patterns but we will also need to consider one-sided patterns at times. A
component of a pattern which is a circle in the interior of a 2-simplex will
be called a{\it \ trivial circle}.

\centerline{\BoxedEPSF{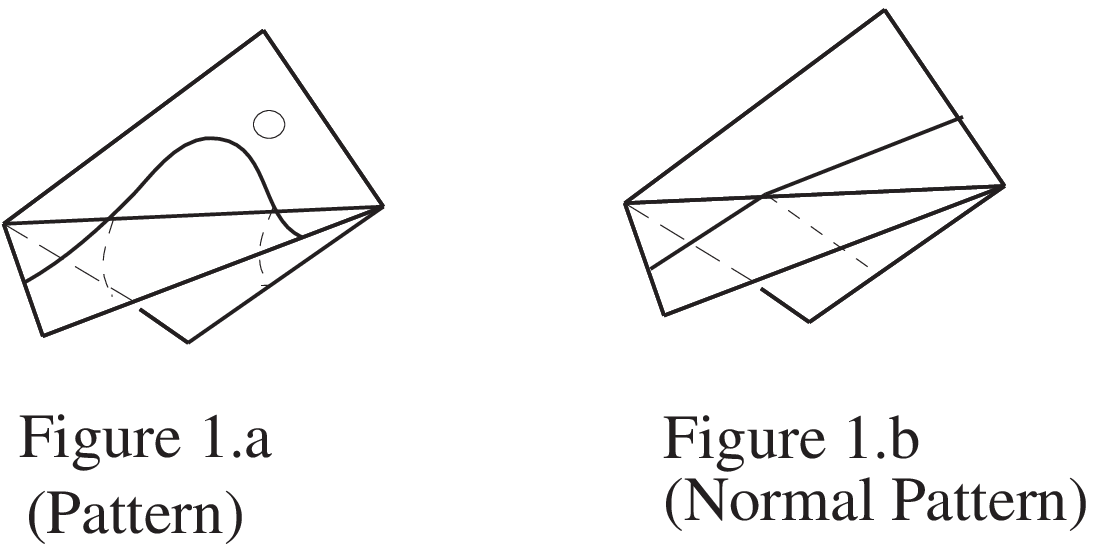 scaled 800}}

Note that a pattern is compact if and only if it meets only a finite number
of the simplices of $Y$. As we are mainly interested in the combinatorics
described by a pattern, we will often say that a compact pattern is finite.
Finite two-sided patterns typically arise from maps to the real line by
taking inverse images of regular points. They can also be constructed using
coboundaries of 0-cochains. A pattern will be called {\it normal} if it
intersects each 2-simplex only in arcs, and the two ends of each arc lie on
different edges of the simplex. See Figure 1b. A connected pattern is called
a {\it track}. It should be noted that our definitions of the terms pattern
and track are more general than Dunwoody's definitions in~\cite{mjd:acc2}.
His patterns and tracks are all normal in our terminology. Note that if $Y\;$%
has a free edge $e,$ which means that $e$ is not a face of any 2-simplex,
then a single point of the interior of $e$ is a normal track in $Y.$ We will
not usually be interested in such tracks, but one does need to be aware of
their existence.

As patterns have similar separation properties to those of a codimension-one
submanifold of a manifold, we will say that a two-sided pattern $t$ in a
2-complex $Y$ {\it bounds} $Z$ if $t$ cuts $Y$ into pieces such that the
union of some of them has closure $Z$, and $Z$ contains $t$ but meets only
one side of each component of $t$. A choice of normal direction for each
component of a two-sided pattern $t$ will be called a {\it transverse
orientation} of $t,$ and we will refer to $t$ with a prescribed transverse
orientation as an {\it oriented pattern}. Our interest in patterns and
tracks in $Y$ is because they give a good way of discussing the number of
ends of $Y.$ They are also a natural generalisation of the idea of a normal
surface in a 3-manifold. For our purposes, we will need to consider patterns
which are not normal, although these will play a small role in our
arguments. Consider a finite pattern $t$ in $Y.$ We will be interested in
whether the closures of the components of $Y-t$ are compact or not. Again,
because it is the combinatorics which are important, a component whose
closure is compact will be called finite. This is equivalent to saying that
the component contains only a finite number of vertices of $Y$. We will say
that a finite pattern $t$ in $Y$ {\it splits} $Y$ if it cuts $Y$ into pieces
at least two of which are infinite. Now consider an oriented finite pattern $%
t.$ If $f$ denotes a proper map of the line into $Y$ or a map of the circle
into $Y$ such that $f$ is transverse to $t,$ we define the sign of an
intersection point of $f$ with $t$ by comparing the orientation of the line
or circle with the transverse orientation of $t.$ (Here we use the
definition that a map $f:X\rightarrow Y$ is {\it proper} if the pre-image of
any compact set is compact.) By summing over all intersection points, we
define the intersection number of $f$ with $t.$ As any proper map of the
line into $Y\;$and any map of the circle can be properly homotoped to be
transverse to $t,$ and as the intersection number obtained is independent of
the homotopy, this defines the intersection number of any such map with $t.$
We say that an oriented finite pattern $t$ is {\it essential} if any loop in 
$Y$ has zero intersection number with $t$ and there is a proper map of the
line into $Y$ which has non-zero intersection number with $t.$ Note that if $%
t$ is not connected, then the condition of essentiality may well depend on
the choice of transverse orientation. We will say that a two-sided finite
pattern $t$ is {\it essential }if $t$ admits an orientation which makes it
essential. Clearly an essential pattern must split $Y.$ The additional
conditions which we impose for a pattern to be essential are needed because
we want to consider patterns which are not connected. They are aimed at
simplifying the discussion of essential patterns but they do not really
restrict us in anyway. To see this, consider a finite pattern $t$ which
splits $Y.$ Next pick an infinite component of the complement of $t$ and
denote its closure by $U.$ We define a sub-pattern $t^{\prime }$ of $t$ to
consist of those components of $t$ which meet $U$ but do not lie in the
interior of $U.$ Then $t^{\prime }$ bounds $U$ and if we transversely orient 
$t^{\prime }$ to point into $U,$ then $t^{\prime }$ is essential in $Y.$
Thus a splitting pattern always contains an essential pattern. We will say
that an essential pattern $t$ is {\it elementary} if $t$ bounds a subset $U$
of $Y$ such that all the orientations point into $U$ or all point out of $U.$
(Of course, such $U$ must be infinite.) We have just seen that any splitting
pattern contains an elementary essential pattern.

We will say that two oriented finite patterns $t_1$ and $t_2$ in $Y$ are 
{\it equivalent} if, for any proper map $h$ of the real line into $Y$ and
any map $h$ of the circle into $Y,$ the intersection number of $h$ with $t_1$
equals the intersection number of $h$ with $t_2.$ If $t_1$ and $t_2$ are
disjoint, they are equivalent if and only if the union of $t_1$ and of $t_2$
with opposite orientation bounds a compact subset $U$ of $Y$ such that for
each component $V$ of $U,$ either all the transverse orientations point into 
$V$ or all point out of $V.$ Thus the relation of equivalence could
reasonably be called oriented cobordism. We will say that an oriented
pattern is {\it trivial} if it is equivalent to the empty set.

Before we discuss ends, we prove the following fundamental lemma about
patterns. An invariant of a finite pattern $t$ which will be used often is
its {\it weight} $w(t)$, which we define to be the total number of points of 
$t\cap Y^{(1)}$, where $Y^{(1)}$ denotes the 1-skeleton of $Y$.

\begin{lem}
Any finite oriented pattern in $Y$ is equivalent to a finite oriented normal
pattern, possibly empty.
\end{lem}

\begin{rem}
If we apply the procedure below to a one-sided pattern, we will obtain a
one-sided pattern, but equivalence is not defined for one-sided patterns.
\end{rem}

{\noindent {\it Proof: }} (See~\cite{sw:access})Let $t$ be a finite pattern
in $Y$. If $t$ is not normal, then there is a 2-simplex $\sigma $ of $Y$ and
a component of $t\cap \sigma $ which is either a circle $C$ or an arc $%
\lambda $ with both endpoints on one edge $e$ of $\sigma $. If there is such
a circle $C$, we alter $t$ by deleting $C$. If there is such an arc $\lambda 
$, let $\lambda ^{\prime }$ denote the sub-arc of $e$ with the same
endpoints as $\lambda $. We alter $t$ by removing a small neighbourhood of $%
\lambda $ and, for each 2-simplex $\sigma ^{\prime },$ other than $\sigma ,$
which contains the edge $e$, we add to $t$ an arc parallel to $\lambda
^{\prime }$. There is a natural choice of transverse orientation for the new
arc which is compatible with the transverse orientation on the rest of $t.$
See Figure 2.

\centerline {\BoxedEPSF{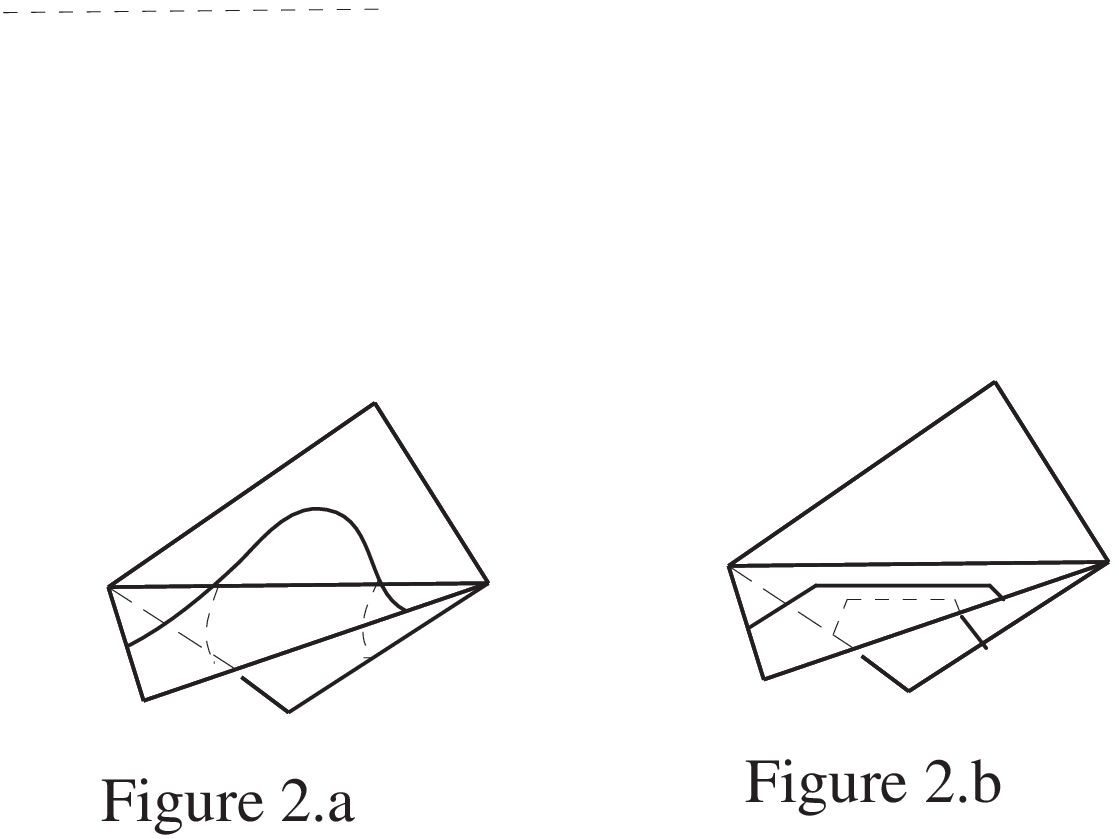 scaled 800}}

Each of these moves replaces $t$ by an equivalent pattern. Also each of
these moves reduces the sum $w(t)+d(t)$, where $w(t)$ is the weight of $t$
and $d(t)$ denotes the sum over all 2-simplices $\sigma $ of $Y$ of the
number of components of $t\cap \sigma $. As this sum is a non-negative
integer, this sequence of alterations to $t$ must terminate, at which point
we will have obtained an oriented normal pattern as required.

Note that even if $t$ has no components which are trivial circles, such
components can be introduced by the second type of alteration. Thus we may
need to delete trivial circles at some stage. Now we consider the connection
between ends and patterns.

\begin{lem}
The locally finite 2-dimensional simplicial complex $Y$ has at least two
ends if and only if it contains a finite essential pattern.
\end{lem}

\begin{rem}
If $Y$ contains a finite essential pattern, the preceding lemma shows that
it contains a finite normal pattern which is essential.
\end{rem}

{\noindent {\it Proof: }}If $Y$ contains a finite essential pattern, it is
immediate that $Y-\{t\}$ has at least two infinite components, so that $Y$
must have at least two ends. Now suppose that $Y$ has at least two ends.
This implies that we can find a proper map $g$ of $Y$ onto the real line $%
{\bf R.}$ Without altering $g$ on the vertices of $Y,$ we can homotop it to
be linear on the simplices of $Y$ and this will be a proper homotopy. Now if 
$z$ is a point of ${\bf R}$ which is not the image of a vertex, then $%
g^{-1}(z)$ is a finite two-sided pattern $t$ which we give the transverse
orientation induced from the orientation of ${\bf R.}$ Any loop in $Y$ has
zero intersection number with $t,$ as any loop in ${\bf R}$ has zero
intersection number with $z.$ Also if we choose a proper map $\lambda $ of
the line into $Y$ such that $g\circ \lambda $ sends the ends of the line to
distinct ends of ${\bf R,}$ then $\lambda $ will have non-zero intersection
number with $t$ as $g\circ \lambda $ has non-zero intersection number with $%
z.$ This shows that $t$ is essential, as required. Note that $t$ will be
normal automatically.

So far, we have not discussed anything really different from the ideas of
Dunwoody. The first new concept we need to introduce is that of a singular
pattern in $Y.$ Note that a pattern $t$ in $Y$ has a natural induced
structure as a 1-complex with vertex set consisting of $t\cap Y^{(1)}.$ Any
circle component of $t$ which lies in the interior of a 2-simplex of $Y$
will not have any vertices, but we will allow this in our definition of a
1-complex for the purposes of discussing singular patterns.

\begin{defn}
A singular pattern in $Y$ is a proper map $f$ of a 1-complex $t$ into $Y$
such that:

\begin{enumerate}
\item  $f$ maps vertices of $t$ to interior points of edges of $Y$,

\item  $f$ maps non-vertex points of $t$ to interior points of 2-simplices
of $Y,$

\item  If $v$ is a vertex of $t$ of valence $d,$ it is mapped to an edge $e$
of $Y$ also of valence $d,$ with distinct rays of the neighborhood of $v$
being mapped to distinct 2-simplices of $Y.$
\end{enumerate}
\end{defn}

We define the {\it weight} $w(f)$ of a singular pattern $f$ to equal the
number of vertices of $t.$ If $f$ is an embedding this agrees with the
original definition of weight. We also define a singular pattern $f$ to be 
{\it normal} if $t$ has no component which is a circle without vertices, and
if, for each edge $\lambda $ of $t,$ the two endpoints of $\lambda $ are
mapped to distinct edges of $Y.$ Again this agrees with the original
definition of normality when $f$ is an embedding.

The conditions in the above definition imply that a singular pattern has
similar local separation properties to those of a pattern. In particular it
has an induced normal bundle, which means that there is a unique bundle over 
$t$ with fibre the interval $I$ such that $f$ extends to the total space of
this bundle and the extended map is injective on each fibre of the bundle.
We will say that $f$ is two-sided if its normal bundle is trivial. In this
case, a choice of normal direction for each component of $t$ will be called
a transverse orientation of $f$ and we will say that $f$ is oriented. Our
main interest in singular patterns is whether a finite singular pattern is
essential in a sense analogous to that for an embedded pattern. We consider
a finite oriented singular pattern $f:t\rightarrow Y.$ We will say that $f$
is {\it essential} if any loop in $Y$ has zero intersection number with $f$
and there is a proper map of the line into $Y$ which has non-zero
intersection number with $f.$ This clearly agrees with the original
definition of essentiality if $f$ is an embedding. Finally we will say that
two singular finite oriented patterns $f$ and $g$ are {\it equivalent} if,
whenever $h$ is a proper map of the real line into $Y$ or a map of the
circle into $Y,$ the intersection number of $h$ with $f$ equals the
intersection number of $h$ with $g.$ Again this agrees with the definition
of equivalence of embedded oriented patterns, when $f$ and $g$ are
embeddings.

Now we want to bring in ideas of Jaco and Rubinstein \cite{jr:pl}. They
defined a PL-analogue of area for a normal surface in a triangulated
3-manifold. We will follow their ideas to define an idea of length for a
finite, possibly singular, pattern in $Y$. Jaco and Rubinstein~\cite{jr:pl}
considered the case when the 2-complex is the 2-skeleton of a 3-manifold but
their ideas extend to the general case. First one needs a ``metric'' on $Y$.
To provide this, we simply choose an identification of each 2-simplex of $Y$
with an ideal triangle in the hyperbolic plane. We then glue these triangles
by isometries to obtain our ``metric'' on $Y$. We will call this a {\it %
hyperbolic structure} on $Y$. Really this structure is defined on $Y-Y^{(0)}$%
, where $Y^{(0)}$ denotes the set of vertices of $Y$. It does yield a metric
on $Y-Y^{(0)}$, but we will not use this fact. Note that there is no
requirement that this metric be complete, and completeness will be
irrelevant for our arguments. For simplicity, we will give the definitions
for embedded patterns only, with brief comments about how they extend to
singular patterns. We define the complexity $c(t)$ of a finite pattern $t$
in $Y$ to be the ordered pair $(w(t),L(t))$, where $L(t)$ denotes the sum
over all 2-simplices $\sigma $ of $Y$ of the length of the 1-manifold $t\cap
\sigma .$ (Note that $L(t)$ is non-zero unless $t$ consists only of isolated
points on free edges of $Y.)$ For a singular pattern $f:t\rightarrow Y,$
note that one does not want to measure the length of the image of $f,$ but
the length with appropriate multiplicity. This is done by pulling back the
metric on $Y$ to a metric on $t.$ It is assumed that we are only considering
finite patterns whenever we talk of complexity. These complexities are
ordered lexicographically and they are added component-wise. Note that a
hyperbolic structure on $Y$ immediately induces a hyperbolic structure on
every cover of $Y.$ The idea now is that if we consider patterns whose
complexity is as small as possible, then such patterns should have good
properties. In particular they should intersect in a reasonably simple way.
This is by analogy with the very nice properties of surfaces of least area
in 3-manifolds \cite{fhs:least} and of shortest curves on surfaces \cite
{fhs:shortest}. It will be very convenient to refer to a pattern of least
possible complexity as {\it shortest}. However, the reader should keep clear
the fact that this term means least complexity $c(t)$ and not least length $%
L(t).$

We say that two patterns $t_1$ and $t_2$ in $Y$ intersect {\it transversely}
if $t_1\cap t_2$ does not meet the 1-skeleton of $Y,$ and for each 2-simplex 
$\sigma $ of $Y$ the 1-manifolds $t_1\cap \sigma $ and $t_2\cap \sigma $
intersect transversely in the interior of $\sigma $. If $t_1$ and $t_2$ are
finite patterns in $Y$ which intersect transversely, we can define a new
pattern $t$ by performing cut and paste at each point of $t_1\cap t_2$ and
rounding corners. As rounding corners reduces length, we see that $c(t)\leq
c(t_1)+c(t_2),$ with strict inequality so long as $t_1\cap t_2$ is not
empty. (For singular patterns which intersect transversely, the same
construction can be made so long as they are in general position. This means
that no point of $f_1(t_1)\cap f_2(t_2)$ is a multiple point of $f_1$ or $%
f_2.)$ This essentially trivial observation will have important consequences
for us in the cases in which we are interested. Note that even if $t_1$ and $%
t_2$ are normal, the new pattern $t$ will usually not be normal. In fact,
for each point $x$ of $t_1\cap t_2,$ exactly one of the two possible cuts
and pastes at $x$ will yield normal arcs, so that there is a unique cut and
paste of $t_1$ and $t_2$ which yields a normal pattern. In Dunwoody's work 
\cite{mjd:acc2}, he used only this cut and paste which explains why he did
not need to consider non-normal patterns. We will definitely need to use
other cut and paste operations. We introduce one further piece of
terminology. If $t_1$ and $t_2$ are finite patterns in $Y$ which intersect
transversely in a non-empty set, then $t_1\cap t_2$ cuts $t_1$ and $t_2$
into pieces which are not patterns. We will call such pieces {\it partial
patterns}. The complexity of a partial pattern is well defined in the same
way as the complexity of a pattern, and if $s_1,\ldots ,s_k$ denote the
partial patterns which are contained in $t_1,$ we have $c(t_1)=\sum c(s_i).$
Note that this equality holds because $t_1\cap t_2$ does not meet any edge
of $Y.$

We want to use cut and paste constructions on shortest patterns to show that
such patterns have nice intersections and we will also need to show that
shortest patterns exist in suitable classes of patterns in order to apply
these ideas. A good example to bear in mind is that of an essential simple
closed curve $C$ on an annulus $A$ equipped with some metric. If this metric
is reasonable, e.g. if $A$ covers a closed surface $F$ and the metric on $A$
is lifted from that on $F$, we can find a shortest loop in the homotopy
class of $C.$ Of course, we do not want to think of the homotopy class of a
pattern. For us a finite pattern is interesting if it is essential. Given a
finite essential pattern $t$ on a 2-complex $Y,$ we will want to know that
there is a shortest pattern equivalent to $t$ and/or that there is a pattern
which is shortest among all essential patterns in $Y.$

Before proving the existence results which we need, we will show how to use
essential patterns of least complexity. The following result will be useful.

\begin{lem}
Suppose that $Y$ has at least two ends, and let $t$ be a finite essential
pattern in $Y$ of least possible weight and without trivial circles. Then $t$
is normal and intersects each 2-simplex of $Y$ in at most one arc.
\end{lem}

{\noindent {\it Proof: }} The fact that $t$ is normal follows from the proof
of Lemma 2.1 as any normalising moves will reduce the weight of $t.$ If $%
E,E^{*}$ are infinite disjoint sets of vertices of $Y$ with $Y^{(0)}=E\cup
E^{*},$ then we can construct from $\delta E$ a pattern which has the
required property. If we have an essential pattern $t$, we take the
corresponding subdivision $E,E^{*}$ of $Y^{(0)}$. Then $\delta E$ gives an
essential pattern $s$ and any 2-simplex which $s$ meets must also be met by $%
t.$ It follows that $w(s)\leq w(t).$ The minimality of $w(t)$ implies that $%
s $ and $t$ have the same weight, and now it follows that $t$ must meet each
2-simplex of $Y$ in at most one arc.

Now we come to a key property of essential patterns of least complexity.

\begin{lem}
Suppose that $Y$ has at least two ends and let $s$ and $t$ be essential
patterns in $Y$ which are shortest among all such patterns. Then either $s$
and $t$ coincide or they are disjoint.
\end{lem}

{\noindent {\it Proof: }}First note that $s$ and $t$ are two-sided and must
have the same complexity. Next recall that each must contain an elementary
essential pattern, so the fact that each is shortest implies that each is
already elementary. Now suppose that $s$ and $t$ are not disjoint and
intersect transversely. Denote the two pieces of $Y$ bounded by $s$ by $S$
and $S^{*},$ and the two pieces of $Y$ bounded by $t$ by $T$ and $T^{*}.$
Thus $s\cap t$ cuts $s$ into two partial patterns $s\cap T$ and $s\cap T^{*}$
which we denote $s^{\prime }$ and $s^{\prime \prime },$ and it cuts $t$ into 
$t\cap S$ and $t\cap S^{*}$ which we denote by $t^{\prime }$ and $t^{\prime
\prime }.$ None of these four partial patterns need be connected even if $s$
and $t$ are connected. Note that $c(s^{\prime })+c(s^{\prime \prime
})=c(s)=c(t)=c(t^{\prime })+c(t^{\prime \prime }).$ Of these four partial
patterns, we consider one of the shortest, i.e. of least complexity. Without
loss of generality, $s^{\prime }$ is a shortest one. Then consider the
patterns $s^{\prime }\cup t^{\prime }$ and $s^{\prime }\cup t^{\prime \prime
}.$ After rounding corners, we obtain two patterns each of which has
complexity strictly less than $c(t),$ which implies that both patterns are
inessential in $Y.$ The pattern $s^{\prime }\cup t^{\prime }$ bounds $S\cap
T $ and its complement, and this complement must be infinite as it contains $%
S^{*}.$ Hence $S\cap T$ must be finite or $s^{\prime }\cup t^{\prime }$
would be essential. Similarly $S^{*}\cap T$ must be finite. But this implies
that $T$ is finite so that $t$ is inessential. This contradiction completes
the proof of the lemma in the case when $s$ and $t$ intersect transversely.

If $s$ and $t$ intersect but not transversely, we use the trick of Meeks and
Yau, see Lemma 1.3 of \cite{fhs:least}, as modified by Jaco and Rubinstein
in section 2.3 of \cite{jr:pl}. Here is how the trick works in our
situation. Note that the intersection of $s$ and of $t$ with each 2-simplex $%
\sigma $ of $Y$ must consist of geodesic arcs in the hyperbolic metric on $%
\sigma .$ We will start by discussing the simplest case when no two geodesic
arcs of $s$ and $t$ coincide, and $s\cap t$ meets the 1-skeleton of $Y$ in
exactly one point $x$ which lies on an edge $e$ of $Y.$ The problem is that
cutting and pasting $s$ and $t$ at $x$ does not yield a strict reduction in
the length of $s\cup t$ as we cannot round corners. We resolve this as
follows. First pick a direction along $e$ and measure the angles at which
the various geodesic arcs of $s$ and of $t$ meet $e$ at $x.$ The fact that $%
s $ is a critical point for the length function implies that the sum of the
cosines of the angles corresponding to the edges of $s$ must be zero, and
the analogous statement holds for the edges of $t.$ In particular, it
follows that $s$ and $t$ must locally cross at $x.$ This means that any
neighbourhood of $x$ contains points of $s$ on each side of $t$ and vice
versa. Now we perturb $s$ slightly to a new track $s_1$ by moving the point $%
s\cap e$ a small distance $\epsilon $ along $e.$ We will discuss how small
to choose $\epsilon $ in a moment. The new track $s_1$ must cross $t$ in at
least one 2-simplex $\sigma _1$ which has $e$ as an edge. Clearly $%
c(s_1)\geq c(s),$ and we can assume the inequality is strict as otherwise $%
s_1$ and $t$ would be two shortest essential tracks on $Y$ which intersect
transversely, contradicting the result of the previous paragraph. Now $%
L(s_1)-L(s)$ is a function of $\epsilon $ which must be zero to the first
order as $s$ is a critical point of $L.$ However, the reduction of length
obtained by cutting and pasting $s_1$ and $t$ in $\sigma _1$ is, to the
first order, linear in $\epsilon $ with (non-zero) coefficient $\mid \cos
\theta -\cos \phi \mid $ where $\theta $ and $\phi $ denote the original
angles in $\sigma _1$ between $e$ and the arcs of $s$ and $t$ which met $e$
at $x.$ It follows that if we choose $\epsilon $ suitably small, then the
reduction will outweigh the increase so that the arguments of the preceding
paragraph will yield a contradiction by applying them to $s_1$ and $t.$ If
no two geodesic arcs of $s$ and $t$ coincide, and $s\cap t$ meets the
1-skeleton of $Y$ in more than one point, essentially the same argument will
work by perturbing $s$ at each of the points where $s\cap t$ meets the
1-skeleton of $Y.$ If some geodesic arcs of $s$ and $t$ coincide, we need to
perturb $s$ along the entire 1-complex $\Gamma $ of coincident arcs and make
the same argument as before at each extreme point of $\Gamma .$ Note that $s$
and $t$ must locally cross at each extreme point $x$ of $\Gamma ,$ so that
there is an intersection between $s_1$ and $t$ in a 2-simplex $\sigma _1$
which contains $x$ and contains non-coincident edges of $s$ and $t$ which
end at $x.$

We will need one further related result for singular essential patterns.

\begin{lem}
Let $f:t\rightarrow Y$ be a finite two-sided essential singular pattern.
Then either $f$ factors through a covering of an essential embedded pattern
in $Y$ or there is an essential embedded pattern $s$ in $Y$ such that $s$
and $f$ are equivalent and $c(s)<c(f).$
\end{lem}

\begin{rem}
This result means that if we can minimise complexity among all singular
patterns in an equivalence class or collection of equivalence classes, then
the least complexity pattern must be an embedding or a covering of an
essential embedded pattern.
\end{rem}

{\noindent {\it Proof: }}We will first consider the case where $f$ is
self-transverse. Then we can perform cut and paste at all the singular
points of $f(t)$ to obtain a pattern $s$ in $Y$ such that $c(s)<c(f).$ In
order to ensure that $s$ is equivalent to $f,$ we choose all the cut and
paste operations to be orientation preserving. This means that we choose the
transverse orientation for $f(t)$ specified by the transverse orientation
which makes $f$ essential and then ensure that each cut and paste preserves
this direction. The orientation preserving condition ensures that any line
or loop which has algebraic intersection number $d$ with $f$ also has
algebraic intersection number $d$ with the new pattern, where we choose the
obvious transverse orientation for the new pattern.

Now we need to consider the more general situation where $f$ need not be
self-transverse. Our aim is to apply the Meeks-Yau trick as in the proof of
the preceding lemma. In other words, we will perturb $f$ by a small amount
to a self-transverse pattern $f^{\prime },$ and then perform cut and paste
on $f^{\prime }$ as in the preceding paragraph to obtain an embedded pattern 
$s$ such that $c(s)<c(f).$ Note that if $f$ is a covering map of an embedded
pattern, there is no way to carry out this procedure. Suppose first that $t$
is connected and that $f$ is not a covering map of an embedded pattern. Then
the methods apply to obtain an embedded pattern $s$ equivalent to $f$ such
that $c(s)<c(f).$ If $t$ is not connected, but $f$ fails to be a covering
map when restricted to some component $t_1$ of $t,$ we can use the same
methods to obtain an embedded pattern $s_1$ equivalent to $f_1,$ the
restriction of $f$ to $t_1,$ such that $c(s_1)<c(f_1).$ Let $F$ denote the
pattern obtained from $f$ by replacing $f_1$ by the inclusion of $s_1$ and
let $\epsilon $ denote $c(f)-c(F).$ We can perturb $F$ on each of the other
components of $t,$ and perform more cut and paste to obtain an embedded
pattern $s$ which is equivalent to $f.$ Further, this can be done while
increasing the complexity of $F$ by less than $\epsilon ,$ so that $%
c(s)<c(f) $ as required.

In the first part of the proof of the preceding lemma, we showed that if $f$
is a two-sided essential self-transverse pattern in $Y,$ then oriented cut
and paste yields an equivalent embedded two-sided pattern. A natural
question is what can one say when one considers a one-sided pattern. The
answer is the following little result which we will use later.

\begin{lem}
Let $f:t\rightarrow Y$ be a finite one-sided pattern in $Y.$ Then either $f$
is an embedding, or there is a one-sided embedded pattern in $Y$ which has
less complexity than $f.$ In particular, $Y$ cannot be simply connected.
\end{lem}

{\noindent {\it Proof: }}If $f$ covers an embedded pattern $s$ in $Y,$ then $%
s$ must also be one-sided. Thus either $f$ is an embedding or $s$ has less
complexity than $f$ as required. Without loss of generality we can assume
that $t$ is connected. Suppose also that $f$ is in general position. Pick a
maximal tree $T$ in $t,$ and let $S$ denote a subset of $t$ consisting of
one point from the interior of each edge of $t-T.$ Then there is a
transverse orientation for $t-S.$ By removing points from $S,$ we can assume
that this transverse orientation does not extend over any point of $S.$ We
can assume that none of the points of $S$ are double points of $f.$ This
allows us to perform oriented cut and paste of $f$ to obtain an embedded
pattern $s$ in $Y$ which is transversely oriented except at the points of $S$%
, and the orientation does not extend across the points of $S.$ It follows
that $s$ is one-sided as required. If $f$ is not in general position, we
apply the Meeks-Yau trick as in the proof of the previous lemma.

Now we want to prove an existence result for patterns of least complexity.
First we consider only patterns of the same combinatorial type as a given
singular normal pattern $f.$ This means that any such pattern is homotopic
to $f$ through singular patterns in $Y,$ and so, in particular, has the same
weight as $t.$ We call such a homotopy a {\it normal} homotopy.

\begin{lem}
Let $Y$ be a locally finite 2-dimensional simplicial complex equipped with
some hyperbolic structure, and let $f:t\rightarrow Y$ be a finite singular
normal pattern in $Y$ such that the restriction of $f$ to each component of $%
t$ cannot be normally homotoped arbitrarily close to a vertex of $Y.$ Then
there is a singular pattern $f^{\prime }$ which is of the same combinatorial
type as $f$ and has least complexity among all such patterns.
\end{lem}

\begin{rem}
The hypothesis on the components of $t$ implies that no component of $t$ can
be a single point.
\end{rem}

{\noindent {\it Proof: }}This is proved exactly as on page 500 of Jaco and
Rubinstein \cite{jr:pl}.

We will be interested in the existence of shortest patterns equivalent to a
given essential singular pattern $f:t\rightarrow Y$ but not necessarily of
the same combinatorial type. In order to apply the above result, we will
need to know that the essentiality of $f$ implies that no component of $t$
can be homotoped arbitrarily close to a vertex of $Y.$ But this need not be
true unless we impose some condition on $Y.$ For example, suppose that $Y$
is obtained from two infinite 2-complexes $Y_1$ and $Y_2$ by glueing them at
a vertex $v$ and let $t_i$ denote the pattern in $Y_i$ which has one vertex
on each edge which contains $v$ and has one edge in each 2-simplex which
contains $v.$ Thus $t_i$ is a copy of the link of $v$ in $Y_i.$ Then each $%
t_i$ is essential in $Y,$ as it separates $Y$ into pieces which are
essentially $Y_1$ and $Y_2,$ and each $t_i$ can be isotoped arbitrarily
close to $v.$

If $Y$ is a locally finite 2-dimensional simplicial complex, we will say
that a vertex $v$ of $Y$ is a {\it splitting vertex} if some component of
the link of $v$ in $Y$ does not bound a compact subset of $Y.$ The
assumption which we will make for the rest of this paper is that the
2-complexes which we consider have no splitting vertices. This will not be a
problem for our applications for the following reason. If $Y$ does have a
splitting vertex, so does its universal cover, and it follows that the
universal cover has at least two ends. Thus we have the following result.

\begin{lem}
Let $Y$ be a locally finite 2-dimensional simplicial complex which covers a
finite simplicial complex $Z$ with fundamental group $G.$ If $Y$ has a
splitting vertex, then $e(G)>1.$
\end{lem}

Now we can state and prove our existence result for shortest patterns.

\begin{lem}
{\it \ } Let $Y$ be a locally finite 2-dimensional simplicial complex
without splitting vertices, and suppose that $Y$ has a hyperbolic structure
lifted from a finite complex $Z$ covered by $Y.$ Let $f$ denote an essential
finite singular pattern in $Y.$ Then:

\begin{enumerate}
\item  There is a singular pattern $f^{\prime }$ which is equivalent to $f$
and has least complexity among all such singular patterns. Further any such $%
f^{\prime }$ must be a covering map of an embedded normal pattern in $Y.$

\item  There is a singular pattern $f^{\prime \prime }$ which has least
complexity among all essential singular patterns in $Y.$ Further, any such $%
f^{\prime \prime }$ must be an embedded normal pattern in $Y.$
\end{enumerate}
\end{lem}

\begin{rem}
In the first part of this lemma, the situation where $f^{\prime }$ covers an
embedded pattern but is not an embedding can certainly occur. Here is a
simple example. Consider any essential track $t$ embedded in $Y$ and let $f$
consist of two parallel copies of $t.$ Then $f^{\prime }$ could map the two
copies of $t$ onto one copy of $t.$
\end{rem}

{\noindent {\it Proof: }} Let $(w_0,L_0)$ denote the infimum of the
complexities of all singular patterns in $Y$ which are equivalent to $f,$
and let $f_i:t_i\rightarrow Y$ denote a sequence of such patterns such that $%
c(f_i)\rightarrow (w_0,L_0)$ as $i$ tends to infinity. As remarked earlier,
the proof of Lemma 2.1 shows that we can assume that each $f_i$ is normal,
as any normalising move reduces complexity. By passing to a subsequence, we
can also assume that each $f_i$ has weight $w_0.$ We consider the
projections of the $f_i$'s into the finite complex $Z.$ This gives us a
sequence $g_i$ of singular patterns in $Z.$ Note that our hypothesis that
the hyperbolic structure on $Y$ is lifted from one on $Z$ implies that the
complexities of $f_i$ and $g_i$ are equal. As there are only finitely many
combinatorial types of singular patterns of a fixed weight in a given finite
2-complex, we can arrange, by again passing to a subsequence, that all the $%
g_i$'s are of one combinatorial type. This means that each $t_i$ can be
identified with a fixed 1-complex $t^{\prime }$ and that for a given vertex $%
v$ of $t^{\prime }$ each $g_i$ maps $v$ to the same edge of $Z.$ Now the
proof of Lemma 2.11 shows that there is a singular pattern $g:t^{\prime
}\rightarrow Z$ of the same combinatorial type as the $g_i$'s and with
complexity $(w_0,L_0)$ unless, for large values of $i,$ there is a component 
$u$ of $t^{\prime }$ whose image $g_i(u)$ lies arbitrarily close to a vertex
of $Z.$ This would imply that $f_i(u)$ also lies very close to a vertex of $%
Y $ and so must cover some component of the link of this vertex. The fact
that $Y$ has no splitting vertices implies that this component bounds a
compact subset of $Y,$ and hence that the restriction of $f_i$ to $t^{\prime
}-u$ is equivalent to $f_i$ which contradicts our assumption that each $f_i$
has the least possible weight $w_0.$ The fact that $g$ has the same
combinatorial type as each $g_i,$ means that $g$ can be lifted to a
(possibly singular) pattern $f_i^{\prime }$ in $Y,$ where $f_i^{\prime }$ is
of the same combinatorial type as $f_i$. For each $i,$ the singular pattern $%
f_i^{\prime }$ is equivalent to the original singular pattern $f,$ and it
has least possible complexity because its complexity is equal to the
complexity of $g$. Now Remark 2.9 shows that $f_i^{\prime }$ must be an
embedding or a covering map of an embedded normal pattern.

For the second part of the lemma, the existence part of the proof is
essentially the same as in the preceding argument, and is carried out by
minimising over the class of all essential singular patterns in $Y.$ Now
consider a singular essential pattern $f^{\prime \prime }:t^{\prime \prime
}\rightarrow Y$ of least possible complexity. We already know that $%
f^{\prime \prime }$ must be a covering map of an embedded normal pattern $s$
in $Y,$ so that $c(s)\leq c(f^{\prime \prime })$ with equality if and only
if $f^{\prime \prime }$ is an embedding. There must be a component of $Y-s$
with infinite closure $U$ and with infinite complement as otherwise $%
f^{\prime \prime }$ could not be essential. Let $s_1$ denote the sub-pattern
of $s$ which is the union of those components of $s$ which meet $U$ but do
not lie in the interior of $U.$ Then $s_1$ is an essential pattern in $Y.$
As $c(s_1)\leq c(s)\leq c(f^{\prime \prime }),$ the fact that $f^{\prime
\prime }$ minimises complexity over all essential patterns in $Y$ implies
that these inequalities are equalities and so $c(s)=c(f^{\prime \prime })$
which implies that $f^{\prime \prime }$ is an embedding as required.

Using the ideas of the preceding proof, one can give a somewhat different
proof of Stallings' Theorem for finitely presented groups and of the
accessibility of finitely presented groups by following Dunwoody's arguments
in \cite{mjd:acc2}. We sketch the proof of Stallings' Theorem. Consider a
finitely presented group $G$ which is the fundamental group of a finite
2-dimensional simplicial complex $Y_G$ with universal cover $Y.$ First we
suppose that $Y$ does not have any splitting vertices. We choose a
hyperbolic structure on $Y_G$ and give $Y$ the induced hyperbolic structure.
The hypothesis of Stallings' Theorem is that $Y$ has at least two ends. Now
we can apply Lemma 2.14 to obtain a shortest essential pattern $t$ in $Y.$ A
crucial point which comes in here is that as $Y$ is simply connected, each
component of $t$ must separate $Y.$ It follows immediately, that some
component of $t$ must be essential, and the fact that $t$ has least
complexity shows that $t$ must equal this component, so that $t$ is
connected. Now we consider the action of $G$ on $Y.$ Lemma 2.7 tells us that
for each $g$ in $G$ the translate $gt$ of $t$ is disjoint from or coincides
with $t.$ It follows that $t$ projects into $Y_G$ as a covering of some
embedded track $s.$ This immediately gives a splitting of $G$ over the
finite group $C$ which is the stabiliser of $t.$ If $Y$ does have a
splitting vertex $v,$ the existence result for shortest tracks fails but it
is not needed. Let $u$ denote the image of $v$ in $Y_G.$ It is immediate
that $u$ determines a splitting of $G$ over the trivial group, which is the
stabiliser of $v.$ Thus in either case, Stallings' Theorem follows for
finitely presented groups.

\section{Ends of pairs of groups}

We will now fix some of the notation to be used in most of the rest of this
paper. $G$ is a torsion free hyperbolic group with one end and $G$ has an
infinite cyclic subgroup $H$ such that $e(G,H)\geq 2$. Following Bowditch,
we sometimes call the number $e(G,H)$ the number of co-ends of $H$ (in $G$). 
$X$ denotes a simply-connected 2-dimensional simplicial complex on which $G$
acts on the left freely and simplicially so that the quotient is compact and
inherits a simplicial complex structure. If $L$ is a subgroup of $G$ then $%
X_L$ denotes the quotient complex $X/L.$ We will always use only the
triangulation of $X_L$ inherited from the given triangulation of $X.$ The
notation $X_H$ usually means that $H$ is infinite cyclic and if $h$ is a
generator of $H$ we sometimes write $X_h$ instead of $X_H$.

We will often be interested in the subgroup of $\pi _1(Y)$ represented by a
track $t$ in $Y,$ where $Y$ is any locally finite 2-dimensional simplicial
complex. Of course, the fundamental group of any track is a free group, but
the group in which we will be interested is the image of the natural map $%
\pi _1(t)\rightarrow \pi _1(Y).$ If $K$ denotes this subgroup, we will say
that $t$ {\it carries }$K$. Note that $K$ is only defined up to conjugacy in 
$\pi _1(Y)$, but this will not cause any problems as $\pi _1(Y)$ will
usually be abelian in the cases in which we are interested.

Now we consider an infinite cyclic subgroup $H$ of $G$ and the corresponding
2-complex $X_H$. We will assume that $H$ has infinite index in $G$, so that $%
X_H$ is not compact. The reader should think of the special case where $X_G$
is a closed orientable surface so that $X_H$ is an annulus which contains an
embedded essential simple closed curve. From our point of view, the
fundamental property of this simple closed curve is that it separates the
annulus into two infinite pieces. In addition, of course, it is connected
and carries the entire group $H.$ In the general situation, we want to find
a track which will play the role of this simple closed curve on the annulus.

\begin{lem}
Let $H$ be an infinite cyclic subgroup of $G$ which has infinite index in $G$%
. If $t$ is a finite two-sided track in $X_H$, then it must separate $X_H$.
If, in addition, $t$ carries the trivial group, then $t$ must separate $X_H$
into two pieces one of which is compact and carries the trivial group, so
that, in particular, $t$ is inessential.
\end{lem}

{\noindent {\it Proof: }} If $t$ does not separate $X_H,$ it follows that $t$
must carry the trivial subgroup of $H$. Hence $t$ must lift to $X$, and as $%
X $ is simply connected this lift $s$ must separate $X$ into two pieces.
Each of these pieces must be infinite, which contradicts our hypothesis that 
$X$ has only one end. See~\cite{sw:ends2} for a related argument. Now we
consider the case when $t$ is separating and carries the trivial group.
Again $t$ must lift to $X$, and this lift $s$ must separate $X$ into two
pieces. If both the components of the complement of $t$ in $X_H$ are
infinite, or if the finite component carries a non-trivial subgroup of $H$,
then it follows that the complement of $s$ in $X$ consists of two infinite
components contradicting our hypothesis that $X$ has only one end. This
completes the proof of the lemma.

Now we will concentrate on finite tracks $t$ in $X_H$ which carry a
non-trivial subgroup of $H.$

\begin{lem}
Let $H$ be an infinite cyclic subgroup of $G$ which has infinite index in $G$%
. Then there is a finite two-sided track in $X_H$ which carries a
non-trivial subgroup of $H$.
\end{lem}

{\noindent {\it Proof: }} As $X_H$ is not compact, we can take a proper map $%
\phi :X_H\rightarrow R^{+}$, the non-negative reals and consider the finite
pattern $t=\phi ^{-1}(c)$ where $\phi ^{-1}([0,c])$ contains a loop
representing $h$. At least one of the components of $t$ must be a track of
the required type, by the preceding lemma.

\begin{rem}
{\em Another way to obtain tracks is to take suitable sets in $X_H$ whose
frontiers are disjoint from vertices.}
\end{rem}

\begin{lem}
If $e(X_H)\geq 2$, there is a finite essential track $t$ in $X_H.$ Any such
track carries a non-trivial subgroup of $H$.
\end{lem}

{\noindent {\it Proof: }} As $e(X_H)\geq 2$, we can take a proper map $%
\varphi $ of $X_H$ onto the reals and consider the finite essential pattern $%
\phi ^{-1}(c)$ for some $c$ in ${\bf R}$. Now we apply Lemma 3.1 to show
that some component of $\phi ^{-1}(c)$ is essential. The same lemma shows
that any essential track in $X_H$ must carry a non-trivial subgroup of $H$.

\begin{rem}
If $W$ is a subset of $X_H$ whose projection to ${\bf R}$ is not onto, we
can choose $t$ to be disjoint from $W$ by choosing $c$ in ${\bf R}$ but not
in the image of $W.$
\end{rem}

The inverse image in $X$ of such a track $t$ in $X$ has one component if $t$
carries $H$. Otherwise $t$ must carry a subgroup of $H$ of finite index and
we have a finite number of components in the inverse image of $t$. If we fix
a generator $h$ of $H$, then we can talk of the positive and negative end
points of any one of these components in $\partial X,$ the Gromov boundary
of $X$ or $G$.

We want to use tracks to subdivide $\partial X$. Note that any track in $X$
is two-sided as $X$ is simply connected. We will generally assume that any
track $A$ in $X$ which we are considering has two ends; often it will have
infinite cyclic stabilizer. The quotient by this stabiliser must be compact
if $A$ has two ends. We will use the term $axis$ for any track in $X$ with
two ends and infinite cyclic stabilizer. We will use $Stab^0A$ to denote the
subgroup of $StabA$ consisting of elements which do not interchange the
sides of an axis $A$, and if $h$ is a generator of $Stab^0A$, then $\partial
A$ or $E(A)$ will denote the fixed points of $h$. Usually $P_A$ denotes the
positive or attracting fixed point of $h$, and $N_A$ denotes the negative or
repelling fixed point of $h$. These will also be called positive and
negative end points of $A$. By Coornaert~\cite{coo:thesis}, $H=Stab^0A$ acts
properly discontinuously on $X-E(A)$ with compact quotient. If $t$ is the
image of $A$ in $X_H$, then $t$ carries $H.$ We will say that the axis $A$
is {\it essential} if $t$ is an essential track in $X_H$. Note that if $A$
is essential, then $e(X_H)\geq 2,$ but the following example shows that the
converse is false.

\begin{ex}
Let $X_H$ be obtained from an open annulus $Z$ and a compact annulus $%
Z^{\prime }$ by identifying one component of $\partial Z^{\prime }$ with an
essential circle $S$ of $Z.$ Consider circles $C$ and $C^{\prime }$ in $Z$
and $Z^{\prime }$ respectively which are parallel to $S$ and disjoint from $%
S.$ By an isotopy, we can arrange that each is a track in $X_H$. Then $C$ is
an essential track in $X_H,$ but $C^{\prime }$ is not essential. Thus the
pre-images of $C$ and $C^{\prime }$ in $X$ are tracks $A$ and $A^{\prime }$
with the same endpoints, with $A$ essential and $A^{\prime }$ inessential.
\end{ex}

If $A$ is an essential axis in $X,$ we denote the two parts of the image of $%
(\partial X-E(A))/H$ separated by $t$ as $\partial _0(X_H)$ and $\partial
_1(X_H)$ and think of these as the left and right parts of the boundary $%
\partial (X_H)$ of $X_H$. If $p$ denotes the projection from $X_H\cup
(\partial X-E(A))/H$ to $X_H\cup \partial (X_H)$,then the subdivision of $%
(\partial X-E(A))$ consists of $p^{-1}(\partial _0(X_H))$ , $p^{-1}(\partial
_1(X_H))$. If we take any finite cover $X_K$ of $X_H$ and the lift $%
t^{\prime }$ of $t$ and repeat the above construction, then we get the same
subdivision of $(\partial X-E(A))$. On the other hand, if we start with a
track $t$ whose fundamental group does not surject onto $H$, then the
subdivision of $(\partial X-E(A))$ may depend on the lift of $t$ that we
take in $X_K.$ Finally note that distinct essential axes $A$ and $A^{\prime
} $ with the same end points may yield different subdivisions of $(\partial
X-E(A)),$ as the following example shows.

\begin{ex}
This example is the same as the preceding one except that we take $Z^{\prime
}$ to be a half open annulus, so that $X_H$ has three ends. Now $C$ and $%
C^{\prime }$ are both essential tracks in $X_H,$ but $A$ and $A^{\prime }$
divide $(\partial X-E(A))$ in distinct ways.
\end{ex}

The positive outcome of the above discussion is the following.

\begin{lem}
Given an essential axis $A$ in $X$ and an infinite cyclic $H\subseteq Stab^0A
$, then $A/H$ defines a subdivision of $\partial X_H$ and the corresponding
subdivision of $(\partial X-E(A))$ is independent of the $H$ chosen in $%
Stab^0A$.
\end{lem}

If $A,A^{\prime }$ are disjoint and $H\subseteq Stab^0A\cap Stab^0A^{\prime
} $ and if $A,A^{\prime }$ bound a compact subset of $X_H$, we get the same
subdivision of $(\partial X-E(A))$ whether we use $A$ or $A^{\prime }$ (note
that $E(A)=E(A^{\prime })$). This generates an equivalence relation on
essential axes with the same end points. We will denote these subsets of $%
(\partial X-E(A))$ by $\partial _0A$ and $\partial _1A$ ; thus implicitly we
are using some transverse orientation for $A$. We will also often use the
notation $\partial _LA$ and $\partial _RA$ for the above sets and $L_A$,$R_A$
for the corresponding parts of $X$ separated by $A$.

\begin{rem}
{\em We will also use another description of the subdivision of $\partial X$%
. Taking a suitable base point $x$ in $R_A$, we can describe $\partial _LA$
as the set of points of $\partial X$ which can be joined by quasi-geodesic
rays from $x$ which cross $A$ an odd number of times and $\partial _RA$ as
the set of points of $\partial X$ which can be joined by quasi-geodesic rays
from $x$ crossing $A$ an even number of times.}
\end{rem}

\section{Crossing of Tracks}

We next consider crossing of tracks in $X$. It is understood that all the
tracks considered are two-ended and when we talk of subdivision of the
boundary, the track considered has infinite cyclic stabilizer.

\begin{defn}
We say that $B$ crosses $A$, if the end points of $B$ are on different sides
of $A$, that is, one end point of $B$ is in $\partial _LA$ and the other in $%
\partial _RA$.
\end{defn}

It is clear that if $B$ crosses $A$ and $B^{\prime }$ is a track with the
same end points as $B,$ then $B^{\prime }$ also crosses $A.$ Also if $A$ is
equivalent to $A^{\prime },$ meaning that $A$ and $A^{\prime }$ induce the
same splitting of $\partial X$, then $B$ and $B^{\prime }$ cross $A^{\prime
} $. It is easy to see that crossing is not symmetric in general even if one
restricts attention to axes (see the examples later in this section). This
will cause us some difficulty, but it turns out that we will have symmetry
in the cases in which we are interested. We note that if $B$ crosses $A$,
then any bi-infinite quasi-geodesic joining the end points of $B$ crosses $A$
an odd number of times. We can, if necessary, modify the quasi-geodesic to
cross $A$ exactly once.

We now examine the symmetry of crossing. The following proposition shows
that we have symmetry in many cases.

\begin{prop}
Suppose that $A,B$ are essential axes so that $e(X_H)\geq 2$ and $%
e(X_{H^{\prime }})\geq 2$, where $H=Stab^0A$ and $H^{\prime }=Stab^0B$. Then 
$A$ crosses $B$ if and only if $B$ crosses $A$.
\end{prop}

{\noindent {\it Proof: }}It will suffice to show that if $B$ does not cross $%
A$, then $A$ does not cross $B.$ If $B$ does not cross $A,$ the picture in $%
X_H$ is as shown in Figure 3. Remark 2.12 tells us that we can choose an
essential track $t^{\prime }$ in $X_H$ disjoint from $\pi (B)$ where $\pi $
denotes the covering projection $X\rightarrow X_H$. This gives a track $%
A^{\prime }$ in $X$ which is disjoint from $B$ with the same end points as $%
A $. Thus $A^{\prime }$ cannot cross $B$. As $A$ and $A^{\prime }$ have the
same end points, it follows that $A$ cannot cross $B$. This completes the
proof of Proposition 4.2.

\centerline {\BoxedEPSF{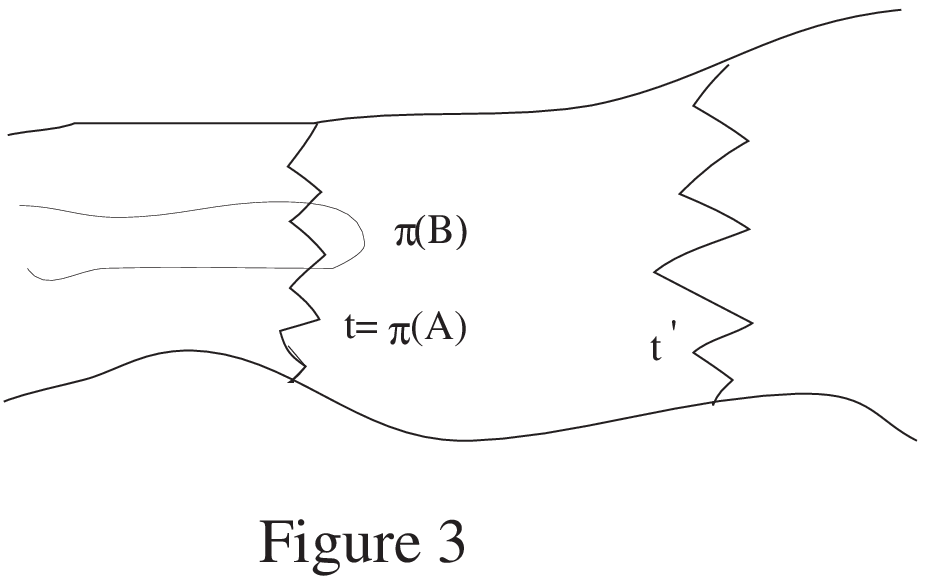 scaled 800}}

The next proposition is crucial to our whole approach.

\begin{prop}
Let $A$ be an essential axis, and suppose that for some $H\subseteq Stab^0A$%
, $e(X_H)\geq 3$. Then the end points of $A$ lie on one side of any axis $B$
in $X$, i.e. $A$ does not cross $B$.
\end{prop}

{\noindent {\it Proof: }} Consider $t=\pi (A)$ in $X_H$. There are three
possibilities for the end points of $\pi (B)$: (1) Both end points are in $%
\partial _0X_H$, (2) One end point is in $\partial _0X_H$ and the other in $%
\partial _1X_H$, (3) Both the end points are in $\partial _1X_H$ (See Figure
4).

\centerline {\BoxedEPSF{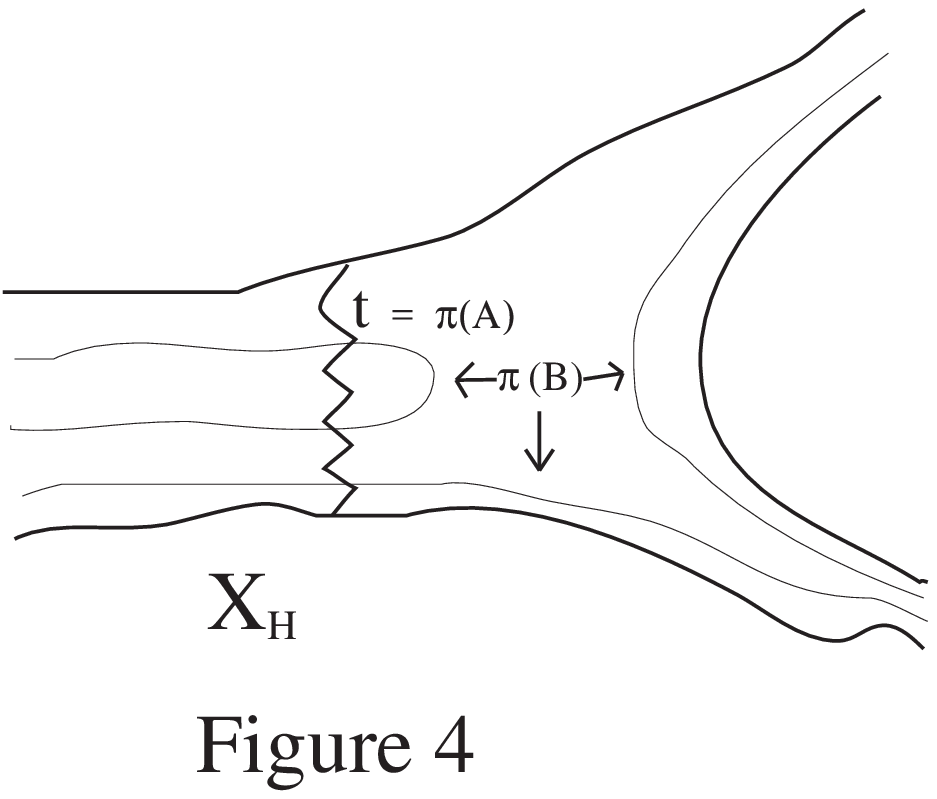 scaled 800}}

In each case, the fact that $X_H$ has at least three ends while $B$ has only
two ends implies that there is a proper map $\phi :X_H\rightarrow {\bf R}$
such that the induced map from $B$ to ${\bf R}$ is not onto. Now Remark 3.5
implies that we can find a finite track $s$ in $X_H$ which carries a non
trivial subgroup of $H$ and does not intersect $\pi (B)$ (see Figures 4a and
4b). This gives a track $A^{\prime }$ in $X$ with the same end points as $A$
which is disjoint from $B$.

\centerline {\BoxedEPSF{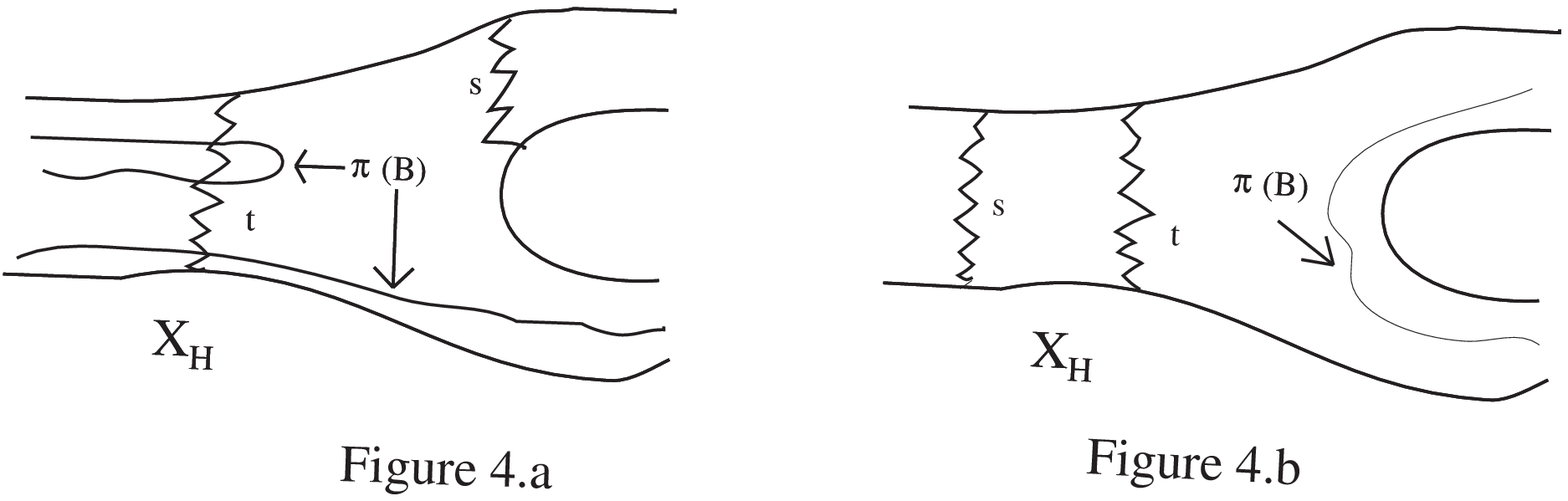 scaled 800}}

It follows that $A^{\prime },$ and hence $A,$ does not cross $B.$ This
completes the proof.

Now we give the promised examples of asymmetry of crossing of axes.

\begin{ex}
{\em For any hyperbolic group }$G${\em \ which is not elementary, if }$A$%
{\em \ is an essential axis with }$H\subseteq Stab^0A${\em , and if }$%
e(X_H)\geq 3,${\em \ there is an axis }$B${\em \ such that }$B${\em \
crosses }$A${\em \ but }$A${\em \ does not cross }$B.${\em \ This is
immediate if we can find }$B${\em \ as in case 2) in the proof of
Proposition 4.3. To find such a }$B${\em , we take }$h\in G${\em \ which has
limit points as in case 2), (such }$h${\em \ must exist by the density of
endpoints of axes \cite{tukia:conv}), and take }$B${\em \ to be the axis of }%
$h${\em . By Proposition 3.2 we must have }$e(X_{h^n})=1${\em \ for any
positive integer }$n${\em . Similar examples can be constructed even when }$%
e(G,H)\leq 2${\em \ for all two-ended }$H${\em \ in }$G${\em \ (See the
definition and examples below).}
\end{ex}

This example may give the impression that asymmetry in the crossing of axes
is a surprising condition which can only occur for rather subtle reasons.
The following simple example shows that asymmetry of crossing is not
surprising and should have been expected.

\begin{ex}
{\em Let }$M${\em \ be a compact orientable hyperbolic 3-manifold with an
incompressible boundary component }$F${\em \ which is not a sphere or torus.
Identify the universal cover of }$F${\em \ with the hyperbolic plane }$H^2.$%
{\em \ Let }$C_1${\em \ and }$C_2${\em \ be essential simple closed curves
on }$F${\em \ such that }$C_1${\em \ and }$C_2${\em \ cannot be homotoped in 
}$F${\em \ to be disjoint. Let }$G${\em \ denote the fundamental group of }$%
M,${\em \ let }$g_1${\em \ and }$g_2${\em \ denote generators of the cyclic
subgroups of }$G${\em \ carried by }$C_1${\em \ and }$C_2${\em \ and let }$%
l_1${\em \ and }$l_2${\em \ denote axes in }$H^2${\em \ for }$g_1${\em \ and 
}$g_2${\em \ which lie above }$C_1${\em \ and }$C_2.${\em \ The fact that }$%
C_1${\em \ and }$C_2${\em \ cannot be homotoped in }$F${\em \ to be disjoint
means that }$l_1${\em \ and }$l_2${\em \ must cross, in the sense of
Definition 4.1.}

{\em Suppose that there is an essential annulus }$S_1${\em \ embedded in }$M$%
{\em \ with }$C_1${\em \ as a boundary component, and let }$S_2${\em \
denote an inessential annulus embedded in }$M${\em \ with }$C_2${\em \ as a
boundary component. This means that }$S_2${\em \ is parallel into }$F.${\em %
\ Fix a triangulation of }$M${\em \ such that }$S_1${\em \ and }$S_2${\em \
are normal surfaces in }$M.${\em \ Let }$X_G${\em \ denote the 2-skeleton of 
}$M${\em , and let }$t_1${\em \ and }$t_2${\em \ denote the intersection of }%
$S_1${\em \ and }$S_2${\em \ with }$X_G.$ {\em \ Thus }$t_1${\em \ and }$t_2$%
{\em \ are tracks in }$X_G${\em \ which contain }$C_1${\em \ and }$C_2${\em %
\ and hence carry the cyclic groups generated by }$g_1${\em \ and }$g_2${\em %
\ respectively. Let }$X${\em \ denote the universal cover of }$X_G,${\em \
and let }$A_1${\em \ and }$A_2${\em \ denote components of the pre-image in }%
$X${\em \ of }$t_1${\em \ and }$t_2,${\em \ so that they contain }$l_1${\em %
\ and }$l_2${\em \ and are axes for }$g_1${\em \ and }$g_2${\em \
respectively. The fact that }$S_2${\em \ is inessential implies that }$%
e(G,A_2)\leq 1${\em \ so that }$A_1${\em \ cannot cross }$A_2.${\em \ The
facts that }$S_1${\em \ is essential and that }$l_2${\em \ crosses }$l_1$%
{\em \ implies that }$A_2${\em \ crosses }$A_1.$
\end{ex}

Recall that a one-ended hyperbolic group $G$ is said to have a multi-band or
to be of multi-band type, if there is a two-ended subgroup $H$ with $%
e(G,H)\geq 3$, and that $G$ is of surface type if $e(G,H)\leq 2$ for every
two-ended subgroup $H$. The proof of the main theorem is broken up into two
cases using this division. We will prove the theorem in the multi-band case
in the next section. The rest of the paper will be devoted to providing a
proof in the surface type case. Here are some examples of both types. The
first three are of multi-band type.

\begin{ex}
{\em Consider two orientable closed surfaces of genus $\geq 2$ and identify
them along a non-trivial simple loop on both. If $G$ is the group so
obtained and $H$ the subgroup corresponding to the simple loop, then $%
e(G,H)=4$ and $G$ is hyperbolic by a theorem of Bestvina and Feighn~\cite
{bf:comb}.}
\end{ex}

\begin{ex}
{\em Let }$G_1,G_2${\em \ be hyperbolic groups, let }$g_1,g_2${\em \ be
elements of }$G_1,G_2${\em \ respectively such that }$g_1${\em \ is
indivisible in }$G_1${\em \ and }$g_2${\em \ has a root of order }$\geq 3$%
{\em \ in }$G_2${\em . Let }$G${\em \ be the group obtained by amalgamating }%
$G_1${\em \ and }$G_2${\em \ along the cyclic subgroups generated by }$g_1$%
{\em and }$g_2${\em \ , and let }$H${\em \ denote this cyclic subgroup of }$G
${\em . Again, it is easy to verify that }$e(G,H)\geq 3${\em \ and that }$G$%
{\em \ is hyperbolic by~\cite{bf:comb}.}
\end{ex}

\begin{ex}
{\em This is a special case of the previous example. We can take $G_2$ to be
infinite cyclic and $g_2$ to be the $n^{th}$ power of a generator of $G_2$
where $n\geq 3$. Of course, we take $g_1$ to be indivisible in $G_1$ and so
that it does not generate the whole group. }
\end{ex}

We next construct surface type examples which are not surface groups.

\begin{ex}
{\em Let $M_1,M_2$ be two orientable hyperbolic 3-manifolds with
incompressible boundary. We further assume that there are no essential
annuli or tori in $M_1,M_2$. Identify the manifolds along non-trivial
embedded annuli $A_1\subset \partial M_1,A_2\subset \partial M_2$ to obtain $%
M$. Let $G$ be the fundamental group of $M$ and let $H$ denote the subgroup
corresponding to $A_1,A_2$. We claim that $G$ is of surface type. We have to
show that $e(G)=1$, and $e(G,K)\leq 2$ for any infinite cyclic $K$ in $G$.
It is easy to check that $e(G)=1$, because }$M${\em \ is irreducible with
incompressible boundary. The claim about $e(G,K)$ is divided into a few
cases. If $X$ is the universal cover of $M$ and $X_K$ the cover
corresponding to $K$, both $X,X_K$ are built from various covers of $M_1,M_2$%
. In both cases we need to consider only covers which are either simply
connected or have infinite cyclic fundamental groups. }

{\em If $N$ is any such cover we want to calculate $H_c^1(N;Z)$ and use the
Mayer-Vietoris sequence to estimate $H_c^1(X_K;Z)$. We also have\newline
$H_c^1(N;Z)\cong H_2(N,\partial N;Z)$. The latter fits into the exact
sequence \newline
\[
0\rightarrow H_2(N,\partial N;Z)\rightarrow H_1(\partial N;Z)\rightarrow
H_1(N;Z)...
\]
\newline
Since each $M_i$ is acylindrical, $H_1(\partial N;Z)\cong 0$ or $Z$ and
since $M_i$ is boundary incompressible the map $H_1(\partial N;Z)\rightarrow
H_1(N;Z)$ is injective. Thus $H_c^1(N;Z)$ is always zero in the cases we
need. If $N,N^{\prime }$ are two such copies in $X_K$ and if they intersect,
they intersect along, say $Y$, which is an annulus or the product of an
interval and a real line. Now the exact sequence \newline
\[
0\rightarrow H_c^0(Y;Z)\rightarrow H_c^1(N\cup N^{\prime };Z)\rightarrow
H_c^1(N;Z)\oplus H_c^1(N^{\prime };Z)\rightarrow H_c^1(Y;Z)
\]
\newline
shows that $e(G,K)=1$ when $K$ is not conjugate to a subgroup of $H$, and
that $e(G,K)=2$ if $K$ is conjugate to a subgroup of $H$.}
\end{ex}

\begin{ex}
{\em We can use a construction similar to the above to obtain a
HNN-extension. This time, we start with $M$ with similar properties to those
of $M_i$ above and identify two non-parallel annuli in the boundary of $M$.
The proof that }$M$ {\em is of surface type is similar.}
\end{ex}

\section{Proof in the multi-band case}

We need the following variation of Theorem 2.2 of~\cite{sc1:newpro} stemming
from the work of Stallings \cite{st:spl} and Dunwoody~\cite{mjd:acc1}. Let $%
K $ be a finite simplicial 2-complex with universal cover $X$, and let $A$
be a track in $X$ with stabiliser $H$. Let $\pi $ denote the projection $%
X\rightarrow X_H$ and let $E,E^{*}$ be the two components of $X-A$.

\begin{thm}
Suppose that

(a) $A/H$ is finite, \newline
(b) Both the components of $X_H-A_H$ are infinite, \newline
(c) For every $g\in G=\pi _1(K)$, at least one of the four sets 
\[
gE\cap E,gE\cap E^{*},gE^{*}\cap E,gE^{*}\cap E^{*}
\]

has finite image in $X_H$, and \newline
(d) If two of the four sets in (c) have finite image in $X_H$, then at least
one of them is empty. \newline
(e) the equation $gE=E^{*}$ does not hold for any element $g$ of $G$.\newline
Then $G$ splits over $H$.
\end{thm}

\begin{rem}
Assumption c) means that the axis $A$ does not cross any of its translates.
Assumption d) means that if $gA$ projects to a finite singular track in $X_H$
then $A$ and $gA$ are disjoint or coincide. Assumption e) means that no
element of $H$ interchanges the sides of $A$ in $X$ so that $A/H$ is a
two-sided track, and assumption b) means that $A/H$ is essential in $X_H.$

Ideally, we would like to prove this theorem by taking $A/H$ to be shortest
in some sense and showing that as $A$ does not cross any of its translates,
it must be disjoint from all its translates. This would immediately imply
the result. However,there are some technical problems with this approach so
we simply use the arguments in \cite{sc1:newpro}.
\end{rem}

The theorem in~\cite{sc1:newpro} is stated in terms of 3-manifolds, but the
proof of Theorem 2.2 in~\cite{sc1:newpro} carries over to a proof of Theorem
5.1 above. (The condition (d), is not explicitly stated in \cite{sc1:newpro}%
, though it is mentioned a few times that it is satisfied in the course of
the proof).

In this paragraph, we discuss the case of a general group $G$ which can have
a multi-band or be of surface type. We first examine the situation when one
has two axes $A,B$ which do not cross. We may assume by choosing appropriate
transverse orientations that $N_A,P_A\in L_B$ and $N_B,P_B\in R_A$. Consider
the infinite components of $X-(A\cup B)$. Let $H=Stab^0A$ and let $K=Stab^0B$%
. Choosing a suitable base point near $N_B$ or $P_B$ in $R_B$ (see Figure 9
in the proof of Lemma 8.1) we see that there is an infinite component $X_1$
of $X-(A\cup B)$ whose closure contains $N_B,P_B$ but does not contain
translates by $K$ of $N_A,P_A$. Similarly there is a component $X_2$ whose
closure contains $N_A,P_A$ but does not contain translates by $H$ of $%
N_B,P_B $. Finally, there is either one or two infinite components $X_3,X_4$
(one of them may be empty) such that their closures contain one of the end
points of $A$ and one of the end points of $B$ (possibly all four). The
closures of the $X_i$'s exhaust all of $\partial X$. Thus $X-\cup X_i$ is
compact. With our notation, we have:

\begin{lem}
If $A,B$ are axes which do not cross, then (with the above choices of
transverse orientations) $L_A\cap L_B$, $R_A\cap R_B$ and $R_A\cap L_B$ are
infinite, and $L_A\cap R_B$ is finite.
\end{lem}

We now assume that $G$ has a multi-band and we want to prove Theorem 1.3 by
using Theorem 5.1. We restate the result.

{\bf Theorem 1.3: }{\em Let }$G${\em \ be a one-ended hyperbolic group and
let }$H_0${\em \ be a two-ended subgroup with }$e(G,H_0)\geq 3${\em . Then }$%
G${\em \ splits over a subgroup }$H${\em \ commensurable with }$H_0${\em .}

{\noindent {\it Proof: }} The hypothesis that $H_0$ be two-ended implies
that $H_0$ contains an infinite cyclic subgroup $H_1$ of finite index. Let $%
g $ denote a generator of $H_1$ and let $P$ and $Q$ denote the fixed points
of $g$ in its action on $\partial X.$ Let $H_2$ denote the stabiliser of $%
\{P,Q\}.$ Then $H_2$ must contain $H_1$ as a subgroup of finite index so
that $H_2$ also has two ends. The point of considering $H_2$ is that it is
the unique maximal two-ended subgroup of $G$ which contains $H_0.$ Let $H$
denote the intersection of all the conjugates of $H_1$ in $H_2,$ so that $H$
is normal in $H_2.$ As $H_1$ is of finite index in $H_2$ and is infinite
cyclic, it follows that $H$ is also of finite index in $H_2$ and infinite
cyclic. Let $X_i$ denote the quotient of $X$ by $H_i.$ We denote the
quotient group $H_2/H$ by $K.$ Thus $K$ is a finite group which acts on $X_H$
by covering translations with quotient $X_2.$ As $X_H$ is a finite cover of $%
X_0,$ and $e(X_0)\geq 3$ by assumption, it follows that $e(X_H)\geq 3.$ We
claim that there is a finite essential track $t$ in $X_H$ which is $K$%
-equivariant. This means that, for all $k$ in $K,$ the track $kt$ is
disjoint from or coincides with $t.$ Let $s$ denote a finite essential track
in $X_H,$ and let $N$ denote a $K$-invariant neighbourhood of the union of
all the translates of $s$ by $K.$ Then any component of the frontier of $N$
in $X_H$ is a $K$-equivariant finite track. If no component is essential in $%
X_H,$ then all but one of the components of $X_H-N$ is finite, which would
contradict the assumption that $s$ is essential. Thus some component of the
frontier of $N$ is the required $K$-equivariant essential finite track $t$
in $X_H.$ The pre-image of $t$ in $X$ consists of a finite number of
disjoint axes. Let $A$ denote one of them. The $K$-equivariance of $t$
implies that $A$ is $H_2$-equivariant.

Now we can verify the hypotheses of Theorem 5.1. Let $h$ be an element of $%
G. $ The fact that $e(X_H)\geq 3$ combined with Proposition 4.3 tells us
that $A,$ and hence also $hA,$ cannot cross any axis. Hence if $hA$ has
different end points from $A$, then both the end points of $hA$ are on one
side of $A$, say in $\partial R_A$, and both the end points of $A$ are on
one side of $hA$, say in $h\partial L_A$. Lemma 5.3 tells us that exactly
one of the four sets of 5.1.(c) is finite, namely $L_A\cap hR_A$. The other
three sets have points near the boundary of $X$ and thus project to
unbounded sets under $\pi $. Hence exactly one of the four sets of 5.1.(c)
has finite image in $X_H.$ If $hA$ has the same endpoints as $A$, then $h$
lies in $H_2$ so that $hA$ is equal to $A$ or disjoint from $A.$ This
implies that the conditions (a-d) of 5.1 are satisfied. Condition (e) can
fail only if $t$ covers a one-sided track $u$ in $X_2.$ If this occurs, we
let $v$ denote the boundary of a regular neighbourhood of $u$ in $X_2$ and
replace $t$ by a parallel track $t^{\prime }$ which covers $v.$ The new
track $t^{\prime }$ satisfies condition (e) and automatically satisfies the
other conditions also. This completes the proof of the theorem.

\section{ Orientation in surface type groups}

For the rest of this paper, we will concentrate on surface type groups. This
time we cannot expect to split along a subgroup commensurable with $H$, as
the example of the usual surface groups shows. Tukia~\cite{tukia:conj}
showed how to find better candidates if the original $H$ does not work. We
could combine this approach and Theorem 4.1 to complete the proof. Instead
we follow the ideas of Freedman,Hass, Scott~\cite{fhs:least} and show that
tracks of least complexity, like least area surfaces, result in the desired
splittings.

We recall that we are considering a hyperbolic group $G$ which is a
one-ended torsion free group with $e(G,H)\leq 2$ for every infinite cyclic $%
H\subset G$. As mentioned above this case is very similar to the usual
surface groups and we plan to adopt some of Tukia's ideas in~\cite
{tukia:conj}. We have the additional difficulty of orientation reversing
elements which we discuss first. We start by observing that tracks have some
nice properties in the case of groups of surface type.

\begin{lem}
If $G$ is of surface type and $e(G,H)=2$, and if $t$ is a finite track in $%
X_H,$ then $t$ must be two-sided. If $t$ is essential in $X_H,$ then $t$
must carry $H$.
\end{lem}

{\noindent {\it Proof: }}If $t$ is one-sided, let $W$ denote the space
obtained from $X_H$ by cutting along $t.$ Then $X_H$ has a double cover
which is the union of two copies of $W$ glued along a double cover of $t.$
This double cover must have four ends which contradicts the hypothesis that $%
G$ is of surface type. Hence $t$ must be two-sided. If $t$ is essential,
then Lemma 3.1 shows that $t$ must carry a subgroup $H^{\prime }$ of $H$ of
some finite index $d.$ Thus $t$ splits $X_H$ into two infinite pieces $C$
and $D$ such that $C$ carries $H$ and $D$ carries $H^{\prime }.$ It follows
that the cover of $X_H$ of degree $d$ consists of a finite cover of $C$
together with $d$ copies of $D$ and so it must have at least $d+1$ ends. The
fact that $G$ is of surface type implies that $d$ equals $1,$ so that $t$
must carry $H$ as claimed.

\begin{lem}
Let $G$ be of surface type and $e(G,<g>)=1$ , but \newline
$e(G,<g^n>)=2$ for some $n$. Then $e(G,<g^{2m}>)=2$ and \newline
$e(G,<g^{2m+1}>)=1$.
\end{lem}

{\noindent {\it Proof: }} Let $X_k$ denote the quotient of $X$ by the cyclic
group generated by $g^k.$ Note that the finite cyclic group ${\bf Z}_k$ acts
on $X_k$ with quotient $X_1.$ If $X_k$ has two ends, the action of ${\bf Z}%
_k $ on $X_k$ induces an action on the two ends, yielding a homomorphism
from ${\bf Z}_k$ to ${\bf Z}_2.$ If this homomorphism were trivial, i.e. if
the action preserved the two ends, it would follow that $X_1$ had two ends
also. This contradiction shows that the homomorphism must be non-trivial, so
that $k$ must be even. If $k$ is odd, we deduce that $X_k$ must have one
end. Hence if there is $n$ such that $X_n$ has two ends, then $n$ is even
and the analysis above shows that $X_2$ must also have two ends. Now it
follows that $X_{2m}$ has at least two ends, for all $m,$ and hence that $%
X_{2m}$ has exactly two ends, as $G$ is of surface type. This proves the
assertions of the lemma.

\begin{defn}
If $G$ is of surface type and if $e(G,<g>)=1$ but \newline
$e(G,<g^n>)=2$ for some $n$, then we call $g$ {\bf orientation reversing }.
\end{defn}

If $g$ is not orientation reversing we call $g$ {\bf orientation preserving}.

In the Lemma above, we can map $C$ onto a Mobius band and take the inverse
image of the middle circle to obtain a finite one-sided pattern inside $C$.
Any one-sided component $t$ of this pattern will be a track in $C$ whose
pre-image in the double cover of $C$ is essential. It follows that $t$
carries $<g>.$ This gives an axis in $X$ stabilized by $<g>$ whose sides are
interchanged by $g$. Conversely, if there is such an axis whose sides are
interchanged by $g$, then $g$ is clearly orientation reversing. Thus:

\begin{lem}
Let $G$ be of surface type. An element $g\in G$ is orientation reversing if
and only if there is an axis $A$ in $X$ stabilized by $<g>$ such that the
sides of $A$ are interchanged by $g$.
\end{lem}

We will call $A$ an axis for $g$ if $g$ stabilises $A$. We have:

\begin{lem}
If $G$ is of surface type and $e(G,<g>)=2,$ then all essential axes for $g$
give the same separation of $\partial X$.
\end{lem}

\section{Shortest tracks for groups of surface type}

It is simpler for many of the assertions that follow to assume that $G$ is
of surface type and even though some of the results of the following
sections are valid more generally, {\bf we will continue to confine
ourselves to surface type groups. }Recall that Lemma 6.1 tells us that if $%
e(X_H)=2$, then any finite track $t$ in $X_H$ is two-sided and that if $t$
is essential in $X_H$ then $t$ must carry $H$. This is very similar to the
situation of an annulus cover of a surface as already discussed. We will use
this analogy and use shortest tracks to obtain splittings in the surface
type case. We will need the following existence result for shortest tracks.

\begin{lem}
Let $G$ be a torsion free hyperbolic group of surface type and let $H$ be an
infinite cyclic subgroup. If $e(X_H)=2$, there is an essential pattern $t$
in $X_H$ which is shortest among all such patterns, and any such pattern is
normal, connected (hence a track), two-sided and carries $H$.

Further, there is an infinite cyclic subgroup $H$ of $G$ and an essential
track $t$ in $X_H$ such that $t$ is shortest among all essential tracks in
all quotients of $X$ by an infinite cyclic subgroup of $G.$
\end{lem}

\begin{rem}
The first part of this lema is analogous to the second part of Lemma 2.14.
There is no point discussing the analogue of the first part of Lemma 2.14,
which is about the existence of a shortest pattern in a given equivalence
class as all essential tracks in $X_H$ are equivalent as $X_H$ has two ends.
\end{rem}

{\noindent {\it Proof: }} The second part of Lemma 2.14 shows that there is
an essential singular pattern in $X_H$ which is shortest among all such
patterns and that any such shortest pattern is a normal embedding. Let $t$
denote such a shortest pattern. As we remarked at the start of this section,
Lemma 6.1 shows that each component of $t$ must separate $X_H.$ It follows
that some component of $t$ is essential. As $t$ is shortest it follows that $%
t$ equals this component so that $t$ is connected. Of course, $t$ must be
two-sided as it is essential. Finally, we apply Lemma 6.1 again to show that 
$t$ must carry $H.$

The last part of Lemma 7.1 is proved in much the same way as the existence
result of Lemma 2.14 by projecting into $X_G$ a sequence of tracks whose
complexity approaches the infimum of all possible complexities. The fact
that the tracks being considered are all essential in some quotient of $X$
by an infinite cyclic subgroup of $G$ ensures that some subsequence
converges to a singular track in $X_G,$ and this will yield a possibly
singular track $f$ in some $X_H$ such that $f$ has the required minimal
complexity and is essential, and covers an embedded track $t.$ Now the fact
that $f$ has minimal complexity implies that $f$ must be an embedding and
that $t$ is the required track.

We will also need an analogous existence result for shortest one-sided
tracks.

\begin{lem}
If $X_H$ contains a one-sided track, then there is a shortest one-sided
track in $X_H$. Further, there is an infinite cyclic subgroup $H$ of $G$ and
a one-sided track $t$ in $X_H$ such that $t$ is shortest among all one-sided
tracks in all quotients of $X$ by an infinite cyclic subgroup of $G.$
\end{lem}

{\noindent {\it Proof: }}This result has essentially the same proof as the
existence results for shortest essential tracks in Lemma 2.14. One takes a
minimising sequence of one-sided tracks in $X_H.$ Remark 2.2 shows that we
can assume that each of these tracks is normal. The hypothesis of
one-sidedness replaces the hypothesis of essentiality, and the key fact
needed is that a one-sided track in a 2-complex cannot be very close to a
vertex. For there is a canonical normal direction for a track close to a
vertex $v,$ in which the normal points towards $v.$ One obtains a possibly
singular one-sided track $f$ in $X_H$ of least possible complexity. Lemma
2.10 shows that $f$ must be an embedding.

The second part of the lemma has essentially the same proof.

We next give analogues of results in Freedman, Hass and Scott~\cite
{fhs:least} and \cite{fhs:shortest}. The proofs are essentially identical.
The first result is the analogue of the fact that two shortest simple closed
curves on an annulus must coincide or be disjoint. It is a special case of
Lemma 2.7, once one remembers that a shortest essential track in $X_H$ is
actually shortest among all essential patterns as proved in Lemma 7.1.

\begin{lem}
Suppose $X_H$ has two ends and that $s$ and $t$ are shortest essential
tracks in $X_H.$ Then either $s$ and $t$ coincide or they are disjoint.
\end{lem}

Our next result is the analogue of the fact that a shortest loop on an
annulus lifts to a shortest loop on any finite cover.

\begin{lem}
Suppose $X_H$ has two ends and $t$ is a shortest essential track in $X_H.$
If $X_K$ is any finite cover of $X_H$ and $t^{\prime }$ denotes the lift of $%
t$ to $X_K,$ then $t^{\prime }$ is a shortest essential track in $X_K$.
\end{lem}

{\noindent {\it Proof: }} Let $p:X_K\rightarrow X_H$ be the covering
projection, $\tau $ be a generator of the deck transformation group and let
the order of the cover be $d$. Suppose there is an essential track $%
s^{\prime }$ in $X_K$ with $c(s^{\prime })<c(t^{\prime })$. Lemma 6.1 shows
that $s^{\prime }$ must carry $K.$ Choose $s^{\prime }$ to be a shortest
such track. Consider the translates of $s^{\prime }$ by the powers of $\tau $%
. The preceding lemma implies that each $\tau ^is^{\prime }$ must coincide
with $s^{\prime }$ or be disjoint from $s^{\prime }.$ It follows that $%
s^{\prime }$ projects into $X_H$ covering its image which we denote by $s.$
Clearly $s$ is an essential track in $X_H,$ and so it must carry $H$ by
Lemma 6.1. It follows that each $\tau ^is^{\prime }$ coincides with $%
s^{\prime },$ so that $s^{\prime }$ covers $s$ with degree $d$ and hence $%
c(s^{\prime })=d.c(s).$ As $c(t^{\prime })=d.c(t),$ it follows that $%
c(s)<c(t).$ This contradiction completes the proof of the lemma.

Before continuing, we briefly consider the orientation reversing case. If $g$
is orientation reversing, consider $p:X_{g^2}\rightarrow X_g$, let $\tau $
denote the non-identity covering translation and let $t$ be a shortest
essential track in $X_{g^2}$. Lemma 7.4 shows that $t$ and $\tau (t)$ are
either disjoint or coincide. In the first case $c(p(t))=c(t)$, and in the
second case $2c(p(t))=c(t)$ and $p(t)$ is one-sided in $X_g$. Also recall
from the discussion at the end of section 5 that if one starts with a
one-sided track $t$ in $X_g,$ then its pre-image in $X_{g^2}$ must be an
essential track. Thus

\begin{lem}
If $c(g)$ is the minimum of $c(t)$ for one-sided tracks in $X_g$, and $c(g^2)
$ is the minimum for essential tracks in $X_{g^2}$, then $2c(g)\geq c(g^2)$.
\end{lem}

Now we consider the analogue of the fact that if one has a shortest loop $C$
on an annulus $M$ then its pre-image in the universal cover $\widetilde{M}$
of $M$ is a length minimising line $l$. This means that any compact interval 
$\lambda $ in $l$ is the shortest path in $\widetilde{M}$ connecting the two
points which form $\partial \lambda $. The proof in \cite{fhs:shortest} uses
the facts that $\lambda $ projects injectively into a finite cover $M_1$ of $%
M$ and that the lift $C_1$ of $C$ into $M_1$ is shortest and so any sub-arc
is also shortest in its homotopy class. We need to replace the concept of
homotopy class for our more general situation.

Let $s$ and $s^{\prime }$ be oriented partial patterns in a 2-complex $Y$
which have the same boundary. Thus $s\cup s^{\prime }$ is naturally a
singular pattern in $Y.$ We will suppose that if we take the given
transverse orientation on $s$ and the opposite one on $s^{\prime },$ this
yields a transverse orientation on $s\cup s^{\prime }.$ We will say that $s$
and $s^{\prime }$ are {\it equivalent} if every loop and proper map of the
line into $Y$ has zero intersection number with $s\cup s^{\prime }$ equipped
with this transverse orientation. Now the following result is clear.

\begin{lem}
Suppose that $X_H$ has two ends and let $t$ be a shortest essential track in 
$X_H.$ If $s$ is an oriented partial pattern contained in $t$ which is
equipped with a transverse orientation induced from a transverse orientation
of $t,$ then $s$ is shortest among all oriented partial patterns in $X_H$
which have the same boundary as $s$ and are equivalent to $s.$
\end{lem}

It follows that if the pre-image of $t$ in $X$ is the axis $A,$ then $A$ is
length minimising in the following sense.

\begin{lem}
Suppose that $X_H$ has two ends and let $t$ be a shortest essential track in 
$X_H.$ Let the pre-image of $t$ in $X$ be the axis $A.$ Then any finite
pattern $s$ contained in $A$ is shortest among all partial patterns in $X$
with the same boundary.
\end{lem}

\begin{rem}
Note that as $X$ is simply connected, Lemma 2.10 shows that any singular
pattern in $X$ is two-sided. Thus the union of two oriented partial patterns
in $X$ with the same boundary is automatically orientable. As $X$ has only
one end, it follows that any two partial patterns in $X$ with the same
boundary are equivalent up to change of orientation so that this lemma does
not need any reference to equivalence.
\end{rem}

{\noindent {\it Proof: }}This follows from the preceding lemma because the
finite pattern $s$ in $A$ projects injectively into $X_{g^n}$ for some $n.$

Now we can prove the analogue of the fact that two length minimising lines
in a plane must coincide, be disjoint or intersect transversely at a single
point.

\begin{lem}
If $A,B$ are minimal axes, then they must coincide or be disjoint, or cross.
If they cross then they intersect in a finite graph such that both $A-(A\cap
B)$ and $B-(A\cap B)$ consist of exactly two infinite components each.
Further, given $\in >0,$ there is an $\in $-isotopy of the track $A/Stab^0(A)
$ which arranges that $A$ and $B$ intersect transversely in finitely many
points and that both $A-(A\cap B)$ and $B-(A\cap B)$ consist of exactly two
infinite components each.
\end{lem}

\begin{rem}
This means that minimal axes in $X$ behave very much like lines in the
plane. Note that the lemma does not assert that $A$ and $B$ intersect
transversely although we have no counterexample.
\end{rem}

{\noindent {\it Proof: }}We start with the case when $A$ and $B$ have the
same end points and hence must have a non-trivial common stabiliser $H$
which preserves the two sides of each. The fact that $A$ and $B$ must
coincide or be disjoint follows immediately by applying Lemma 7.4 to the
quotient tracks $A/H$ and $B/H$ in $X_H.$ For the rest of this proof we will
assume that $A$ and $B$ do not have the same end points.

Suppose that $A$ and $B$ are distinct but intersect transversely. It is
automatic that they intersect in a finite number of points. We will consider
the possibility of non-transverse intersection later. If the lemma fails to
hold, there must be a bounded component $R$ of $X-(A\cup B).$ This is true
whether $A$ and $B$ cross or not. Such a region will be bounded by compact
pieces of $A$ and $B$ which will be partial patterns. In the case of lines
intersecting in a plane, one can find $R$ such that its boundary meets each
of $A$ and $B$ in a connected set, but there seems no reason why this should
be possible in general. However, we can still make what is in essence the
usual cut and paste argument. We will perform cut and paste at each point of 
$A\cap B$ which lies in the boundary of $R.$ At each such point one locally
sees four regions exactly one of which is part of $R.$ We choose the cut and
paste which connects this region to the opposite one. We know that cut and
paste operations reduce complexity but this does not mean much as both $A$
and $B$ have infinite complexity. We get round this problem as follows. Let $%
A_{+}$ and $A_{-}$ denote the two infinite pieces of $A-(A\cap B),$ and
similarly for $B_{+}$ and $B_{-}.$ Let $A_1$ denote the partial pattern
obtained by truncating $A$ in $A_{+}$ and $A_{-}$ and describe $B_1$
similarly. Lemma 7.8 tells us that $A_1$ and $B_1$ are each the shortest
patterns in $X$ with their boundary. The result of our cut and pastes is to
replace $A_1$ and $B_1$ by $A_1^{\prime }$ and $B_1^{\prime },$ where $%
A_1^{\prime }$ is obtained from $A_1$ by removing $A\cap R$ and replacing it
with $B\cap R,$ and $B_1^{\prime }$ is obtained similarly. Now we can say
that $c(A_1^{\prime })+c(B_1^{\prime })<c(A_1)+c(B_1).$ This contradicts the
fact that $A_1$ and $B_1$ are each the shortest partial patterns in $X$ with
their boundary, completing the proof of the lemma in the case when $A$ and $%
B $ intersect transversely.

If $A$ and $B$ have non-transverse intersection, we want to apply the
Meeks-Yau trick as in the proof of Lemmas 2.7 and 2.8. If $X-(A\cup B)$ has
a bounded component, we can use this trick and the arguments of the
preceding paragraph to obtain a contradiction. Thus $X-(A\cup B)$ cannot
have a bounded component. Recall from the discussion in the proof of Lemma
2.7 that each extreme point of the finite graph $A\cap B$ is a local
crossing point of $A$ and $B.$ This means that any neighbourhood of such a
point contains points of $A$ on each side of $B$ and vice versa. If $A$ and $%
B$ do not cross, but do intersect, the local crossing property implies that $%
A$ has points on both sides of $B$ which in turn implies that $X-(A\cup B)$
has a finite component, which we know cannot occur. We conclude that if $A$
and $B$ do not cross, then they must be disjoint. If $A$ and $B$ do cross,
the fact that there cannot be a finite component of $X-(A\cup B)$ implies
that each of $A-(A\cap B)$ and $\;B-(A\cap B)$ cannot have a finite
component, so that each consists of exactly two infinite components as
required.

For the last part of the statement of the lemma, we consider the finite
two-sided track $A/Stab^0(A)$ in $X/Stab^0(A).$ Denote the track by $t_1$
and the quotient of $X$ by $X_1.$ As $t_1$ is two-sided in $X_1,$ it has a
neighbourhood homeomorphic to $t_1\times I.$ This defines an isotopy of $t_1$
which moves it to a parallel copy, and all the nearby parallel copies will
automatically be transverse to the image of $B$ in $X_1,$ so that the new
version of $A$ will intersect $B$ transversely in a finite set. Applying the
Meeks-Yau trick again implies that $A\cap B$ must separate each of $A$ and $%
B $ into exactly two infinite pieces as required, so long as we move $t_1$ a
small enough distance.

\section{Subdivisions of the boundary of $X$}

We consider again subdivisions of the boundary, this time with respect to
various crossing essential axes in the surface type case. We have already
noted that if $A$ is an essential axis in $X$, then $\partial L_A,\partial
R_A$ can be described in terms of quasi-geodesic rays from some base point $%
x $. For example, in Figure 9 below,

\centerline {\BoxedEPSF{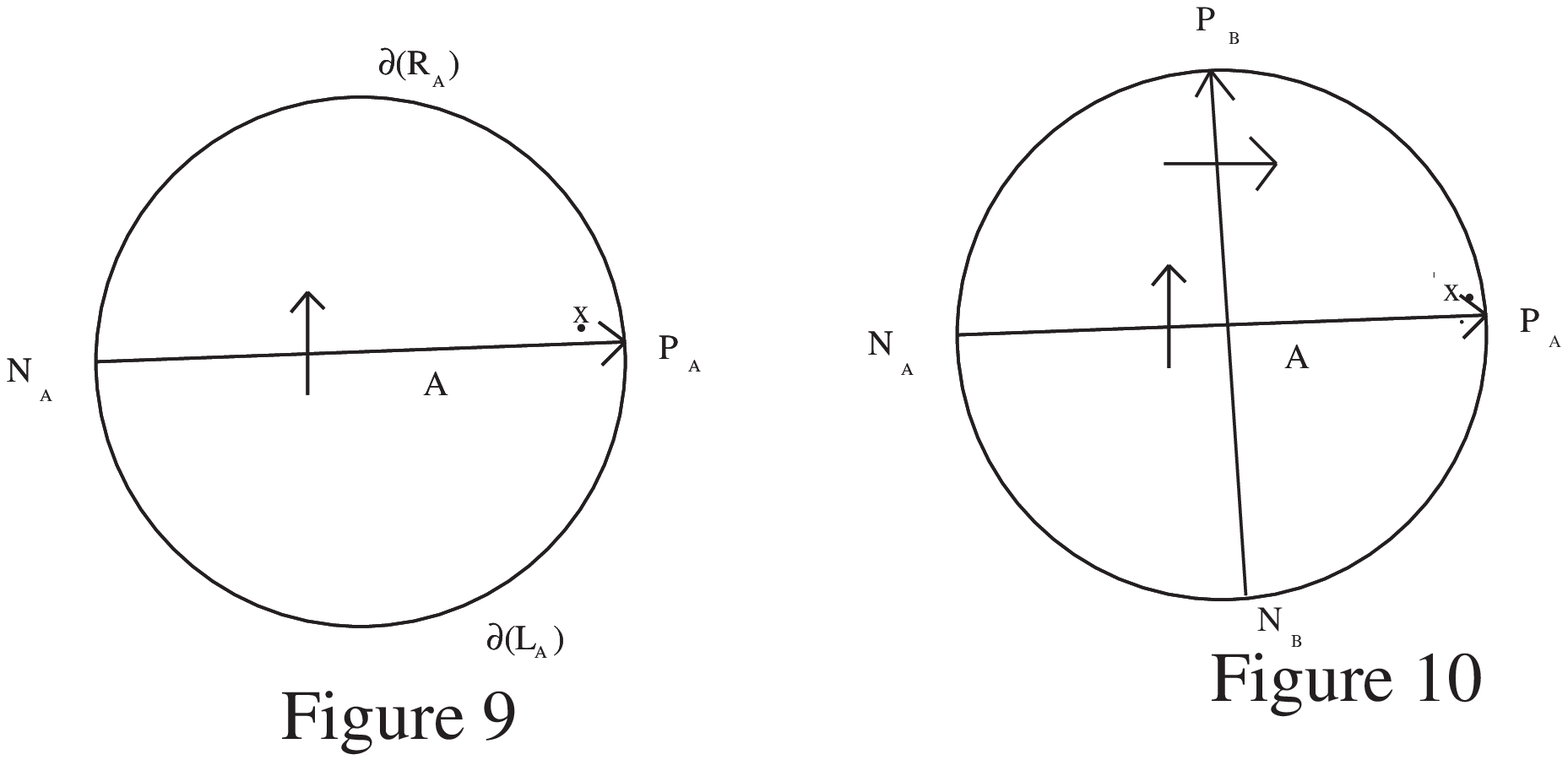 scaled 800}}

$\partial L_A$ consists of points of the boundary for which the
quasi-geodesic rays (hereafter called rays) joining them to $x$ cross $A$ an
odd number of times and in fact we can even find rays which cross exactly
once since $A$ is connected. Similarly $\partial R_A$ is described as the
part of the boundary for which the rays cross an even number of times. We
will use the convention that an arrow representing a transverse orientation
of $A$ always points towards to $R_A.$ Since $X$ has one end we can require
the rays to avoid compact sets in $X$. Suppose $C$ is a compact set in $X$
such that $A-(A\cap C)$ has two infinite components, (there are
automatically at most two). Let these be $E_1,E_2$ with $N_A\in E_1$ and $%
P_A\in E_2$. If $R$ is a ray which avoids $C$ and joins $x$ to a point of $%
\partial L_A$, then the intersection of $R$ with $A$ is in $E_1\cup E_2$. We
can now further modify $R$ to obtain a ray which intersects $A$ exactly
once. We next observe that $L_A,R_A$ have exactly one end each. Thus :

\begin{lem}
We can join $x$ to a point of $\partial R_A$ by a ray which avoids any given
compact set in $X$ and does not intersect $A$. Similarly we can join $x$ to
a point of $\partial L_A$, by a ray which avoids a given compact set $C$ and
meets $A$ exactly once either in $E_1$ or $E_2$ as required.
\end{lem}

Next suppose that two axes $A$ and $B$ cross. We can assume that $P_B$ lies
in $\partial R_A$ by changing our choice of generator of $Stab^0B.$ Figure
10 shows the picture for one choice of transverse orientation of $B$. Given
this choice, we subdivide $\partial X$ into four sets $%
[P_A,P_B],[P_B,N_A],[N_A,N_B],[N_B,P_A]$ as follows: $[P_A,P_B]=\partial
R_A\cap \partial R_B$, $(P_A,P_B)=(\partial R_A\cap \partial
R_B)-\{P_A,P_B\} $ etc. Note that if we change the transverse orientation on 
$B$, then $[P_A,P_B]=\partial R_A\cap \partial L_B.$ Choose a base point $x$
in $X$ which lies in $R_A$ close to $P_A$ and does not lie on $B$. Then $%
(P_A,P_B)$ can be described as the set of points $y$ in $\partial X$ which
can be joined to $x$ by rays which cross both $A$ and $B$ an even number of
times and $[P_A,P_B]$ is its closure. Now join $x$ to $y$ by a ray $r$ which
does not intersect $A,$ and so lies in $R_A$. Such a ray, if it intersects $%
B $ at all, will intersect $B$ in the infinite component of $B-(A\cap B)$ in 
$R_A$ an even number of times and thus can be replaced by a ray which does
not intersect either $A$ or $B$. If $y_1,y_2$ are two points of $(P_A,P_B)$,
then the union of two such rays $r_1,r_2$ is a bi-infinite path joining $%
y_1,y_2$ which does not intersect $A\cup B$. Thus $(P_A,P_B)$ is in the
closure of one infinite component of $X-(A\cup B)$. Hence :

\begin{lem}
If $A,B$ cross then $X-(A\cup B)$ has exactly four infinite components and
their closures intersect $\partial X$ in the subsets 
\[
[P_A,P_B],[P_B,N_A],[N_A,N_B],[N_B,P_A].
\]
\end{lem}

Following Tukia, we sometimes denote these subsets of $\partial X$ by $%
I_1,I_2,I_3,I_4$ in that order. See Figure 10. We next consider the
situation when there is an axis $C$ with $P_C\in (N_A,P_B)$, $N_C\in
(N_B,P_A)$ and observe that we obtain subdivisions of $\partial X$ as in the
case of a circle (See Figure 11 with transverse orientations as shown).

\centerline {\BoxedEPSF{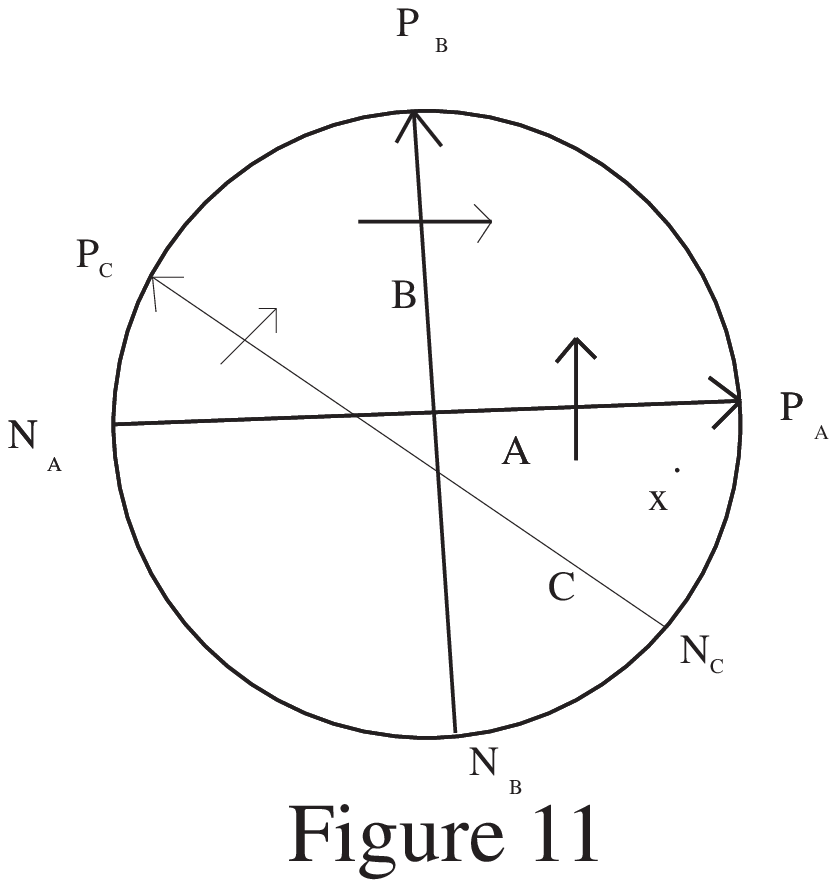 scaled 800}}

We have:

\begin{lem}
Suppose that $C$ is an axis with $P_C\in (N_A,P_B)$, $N_C\in (N_B,P_A)$.
Then we have with suitable choice of transverse orientation for $C$:\newline
(1) $\partial R_C\supset [P_A,P_B],\partial L_C\supset [N_A,N_B]$, \newline
(2) $[N_B,P_A]=[N_C,P_A]\cup [N_B,N_C],[N_C,P_A]\cap [N_B,N_C]=\{N_C\}$, 
\newline
(3) $[N_A,P_B]=[N_A,P_C]\cup [P_C,P_B],[N_A,P_C]\cap [P_C,P_B]=\{P_C\}$ 
\newline
\end{lem}

{\noindent {\it Proof: }} We now take a base point $x$ near $P_A$ in $L_A$
and a ray $R$ joining $x$ to a point of $(P_A,P_B)$. We assume that $R$
avoids all the intersections of $A,B,C$ and does not meet $B$ and meets $A$
exactly once near $P_A$. If $R$ meets $C$ an odd number of times, we modify $%
R$ near the ends of $C$ to $R^{\prime }$ so that $R^{\prime }$ meets $C$
near $P_C$ or $N_C$ and still intersects $A$ once and does not meet $B$.
Since $R^{\prime }$ crosses $A$, $R^{\prime }$ must meet $C$ near $N_C$.
This means that $N_C$ is in $R_A$ and it is in $R_B$ since $R^{\prime }$
does not meet $B$. Hence $N_C$ should be in $R_A\cap R_B$ which contradicts
the hypothesis that $N_C$ is in $(N_B,P_A)$. This proves that $[P_A,P_B]$ is
on one side of $C$. Similarly we conclude that $[N_A,N_B]$ is also on one
side of $C$. We now choose transverse orientation on $C$ so that $%
[N_A,N_B]\subset \partial L_C$ and $[P_A,P_B]\subset \partial R_C$. Since $%
\partial L_C\cap [P_A,P_B]=\emptyset $, we have 
\begin{eqnarray*}
\partial L_C &=&(\partial L_C\cap [N_A,N_B])\cup (\partial L_C\cap
(N_A,P_B))\cup (\partial L_C\cap (N_B,P_A)) \\
&=&[N_A,N_B]\cup (\partial L_C\cap (N_A,P_B))\cup (\partial L_C\cap
(N_B,P_A)).
\end{eqnarray*}
Since $(N_A,P_B)$ is in $\partial L_B,$ we have 
\[
\partial L_C\cap \partial R_B=(\partial R_B\cap [N_A,N_B])\cup (\partial
L_C\cap (N_B,P_A))=(\partial L_C\cap [N_B,P_A)). 
\]
. Hence $[N_B,N_C]\subset [N_B,P_A)$. Since $(P_A,P_B)\subset \partial R_C$
, we have $[N_A,P_C]=\partial L_C\cap \partial R_A=\partial L_C\cap
[N_A,P_B] $. Thus $[N_A,P_C]\subset [N_A,P_B]$. The other inclusions of
Lemma 8.3 are proved similarly.

We note a consequence of the above lemma. Let $g\in Stab^0A$ so that $P_A$
is the attractive fixed point for $g$. We must have the end points of $gB$
in $\partial (R_B)$. For otherwise, let $C=gB$. If $C$ is as in Lemma 8.3,
we see that $g[N_A,P_B]=[N_A,P_C]$ so that $g^i(P_B)\in [N_A,P_B]$ for all
positive $i$ contradicting the assumption that $P_A$ is the positive fixed
point of $g$. If $C=gB$ crosses $B$ in the other direction, i.e. $P_C\in
(P_A,P_B)$, $N_C\in (N_A,N_B)$. we can interchange the roles of $B$ and $C,$
and again arrive at a contradiction. Thus:

\begin{cor}
If $A,B$ are essential axes which cross and if $g\in Stab^0A$, then the end
points of $gB$ are in $\partial (R_B)$, the half containing $P_A$, the
positive fixed point of $g$.
\end{cor}

Up to this point, we have discussed only one choice of configuration for
crossing axes $A$ and $B$ as shown in Figures 10 and 11. If $B$ equals $hA,$
and we choose a transverse orientation for $A$ and the induced one for $B$,
then there are four possible configurations as shown in Figures 12a, 12b,
12c and 12d.

\centerline {\BoxedEPSF{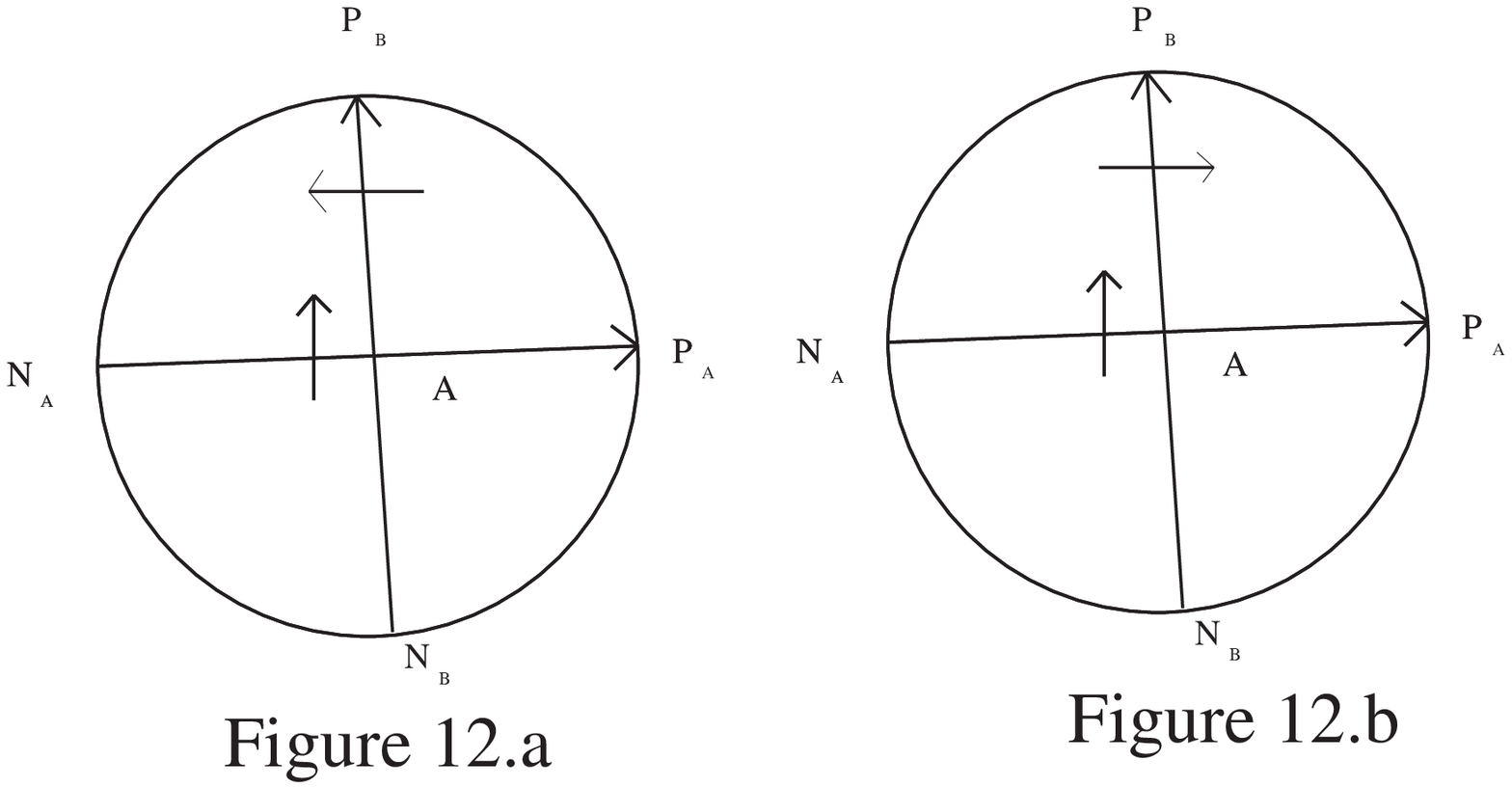 scaled 800}}

\centerline {\BoxedEPSF{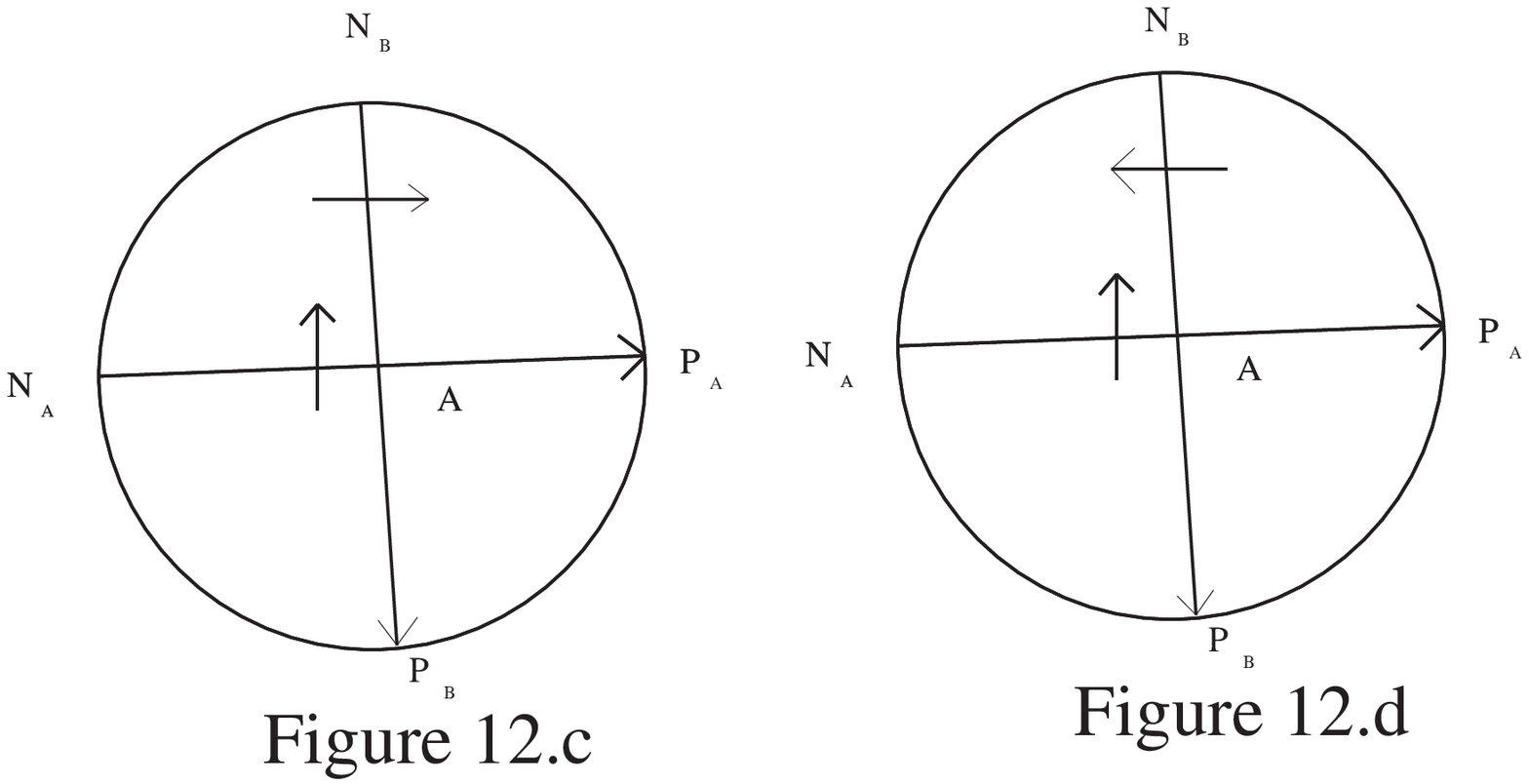 scaled 800}}

Note that configurations 12c and 12d can be obtained from 12a and 12b by
interchanging $A$ and $B$ and replacing $h$ by $h^{-1},$ so that we need
only consider the configurations in 12a and 12b. We will say that $hA$
crosses $A$ in an {\em orientation preserving fashion } in case 12a and in
an {\em orientation reversing fashion } in case 12b. More formally, in the
first case $\{P_B,N_A\}$ ( or $\{P_A,N_B\}$ in case 12c) is on the right or
positive side of both $A$ and $B$ and in the second case $\{P_B,P_A\}$ (or $%
\{N_B,N_A\}$ in case12d) is on the positive side of both $A$ and $B$. We
will use the abbreviations o.r. and o.p. We want to show that if the
crossing is o.p. then $h$ is o.p. Let $A_h$ denote an axis for $h.$ First
observe that the fixed points of $h$ must be in $I_2\cup I_4$ as each must
lie on the positive side of both $A$ and $hA$ or on the negative side of
both. If $A_h$ is not an essential axis, then $h$ is o.p. by definition. So
we may assume that $A_h$ is essential. If both the fixed points of $h$ lie
in $I_2$ or if both lie in $I_4,$ then the endpoints of $A$ and $B$ lie on
one side of $A_h,$ so that $h$ does not interchange the sides of $A_h$ and $%
h $ must be o.p. If one of the fixed points is in $I_2$ and the other in $%
I_4 $, then $A_h$ crosses $A$. Since $A_h$ is an essential axis, crossing is
symmetric so that $A$ crosses $A_h$. If $h$ is o.p., the above corollary
shows that $B=hA$ cannot cross $A,$ a contradiction. If $h$ is o.r., we also
have a contradiction as $P_A$ and $P_B=hP_A$ lie on the same side of $A_h.$
Thus:

\begin{lem}
If $hA$ crosses $A$ in an o.p. fashion, then $h$ is o.p. and the fixed
points of $h$ are in $[N_A,P_B]\cup [N_B,P_A]$. If $h$ has an essential
axis, then the fixed points of $h$ are in either $[N_A,P_B]$ or in $[N_B,P_A]
$ ( that is, $I_2$ or $I_4$).
\end{lem}

We will prove later on (the main difficulty lies in showing that $h$ has an
essential axis) that if the crossing is o.r. then $h$ is o.r. For the moment
we will prove the following lemma.

\begin{lem}
If $hA$ crosses $A$ in an o.r. fashion, then the fixed points of $h$ are in $%
I_1\cup I_3,$ one in $I_1$ and the other in $I_3$.
\end{lem}

{\noindent {\it Proof: }}The fixed points of $h$ must lie in $I_1\cup I_3,$
because if one of them is on the positive (negative) side of $A$ it must lie
on the negative (positive) side of $hA.$ We consider the action of $h$ on $%
\partial X.$ The sets $[N_A,N_B]$ and $[N_B,N_C]$ must be one contained in
the other$.$ Thus $h$ or $h^{-1}$ sends one of these sets into itself and so
this set must contain one of the fixed point of $h.$ In particular $h$ has a
fixed point in $I_3.$ Similarly $h$ must have a fixed point in $I_1.$

\centerline {\BoxedEPSF{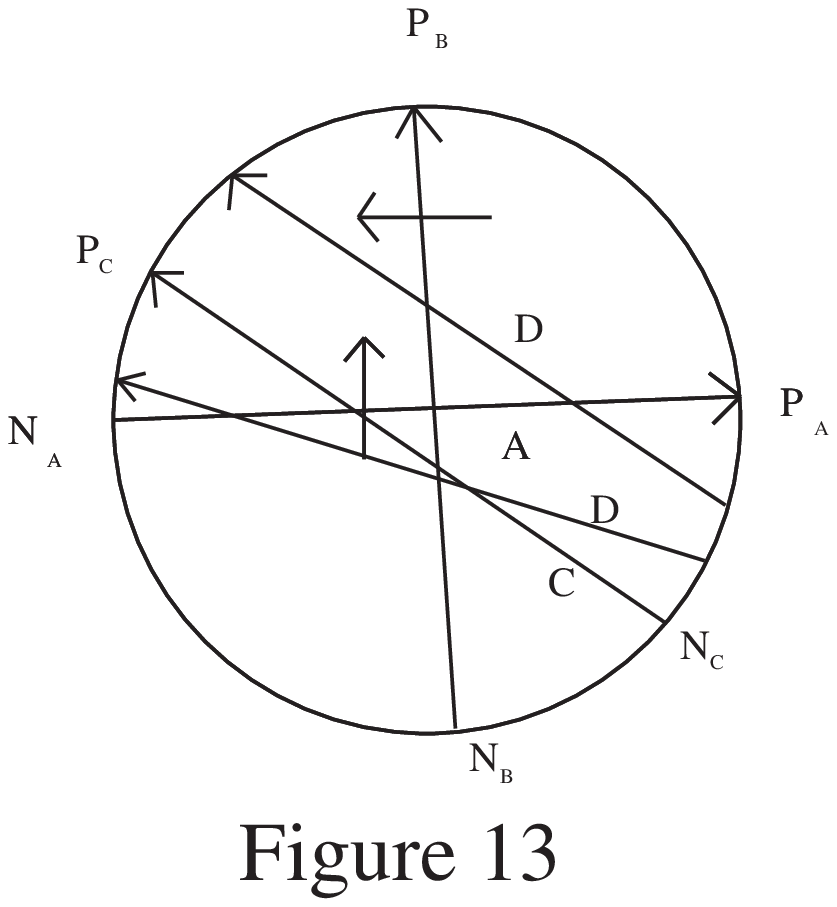 scaled 800}}

Now we consider another axis $D$ which, like $C$ in Lemma 8.3, has $N_D\in
(N_B,P_A)$ and $P_D\in (N_A,P_B).$ See Figure 13. Then Lemma 8.3 tells us
that 
\[
\lbrack N_B,P_A]=[N_B,N_C]\cup [N_C,P_A]=[N_B,N_D]\cup [N_D,P_A]. 
\]
If $N_D\in (N_C,P_A),$ then $N_D$ is above $C$ as in Figure 13. If both
endpoints of $D$ lie above $C$ then clearly $[N_B,N_C]\subset [N_B,N_D]$. In
this case we use the notation $[N_C,N_D]$ for the closure of $%
[N_B,N_D]-[N_B,N_C]$. If $D$ crosses $C$, then replacing $A$ by $D$ in Lemma
8.3, we see that again we have $[N_B,N_C]\subset [N_B,N_D]$. Similar
statements hold if we reverse the roles of $P_D$ and $N_D.$ Thus:

\begin{lem}
If $D_i,1\leq i\leq n$ are a finite number of axes with one end point in $%
[N_B,P_A]$ and the other in $[N_A,P_B]$, then there is a $D_i$ nearest to $%
P_A$ in the sense that $[V_i,P_A]$ does not contain any of the end points of
the other $D_i$'s, where $V_i$ denotes the end point of $D_i$ which lies in $%
[N_B,P_A].$
\end{lem}

\begin{rem}
{\em We apply the above lemma in the following situation. If $B=gA$ and $A,B$
cross, we consider translates $D$ of $A$ with one end point in $(N_B,P_A)$
and the other in $(N_A,P_B)$. Freden~\cite{ef:conv} and Tukia~\cite
{tukia:conv} proved that $G$ acts as a convergence group on $\partial X$. It
follows as in Tukia~\cite{tukia:conj} that there are only finitely many such 
$D$.}
\end{rem}

Following Tukia, we call any translate of $A$ an $A$-axis. We also introduce
the notion of canonical triple. We look at the situation when $A$ and an $A$%
-axis $B$ cross. Suppose that $B=gA$ and that the configuration is as in
Figures 12a or 12b. Then $C=gB$ and $B$ cross. We can vary $g$ by an element
in $Stab^0B$ so that $B$ is unchanged but $C=gB$ is moved up and down. As $%
P_B$ lies in $\partial _RA,$ it follows that $P_C$ lies in $\partial _RB$.
Whether $C$ crosses $A$ or not we can talk of the region $[N_A,P_C]$ in the
first case and $[N_A,N_C]$ in the second case after moving the appropriate
end point of $C$ above $A$. We call the triple $(A,B=gA,C=gB)$ a {\bf %
canonical triple} if $[N_A,P_C]$ ( or $[N_A,N_C]$ in the o.r. case) does not
contain the end point of any $hC$ for $h\in Stab^0B$. See Figures 14a and
14b.

\centerline {\BoxedEPSF{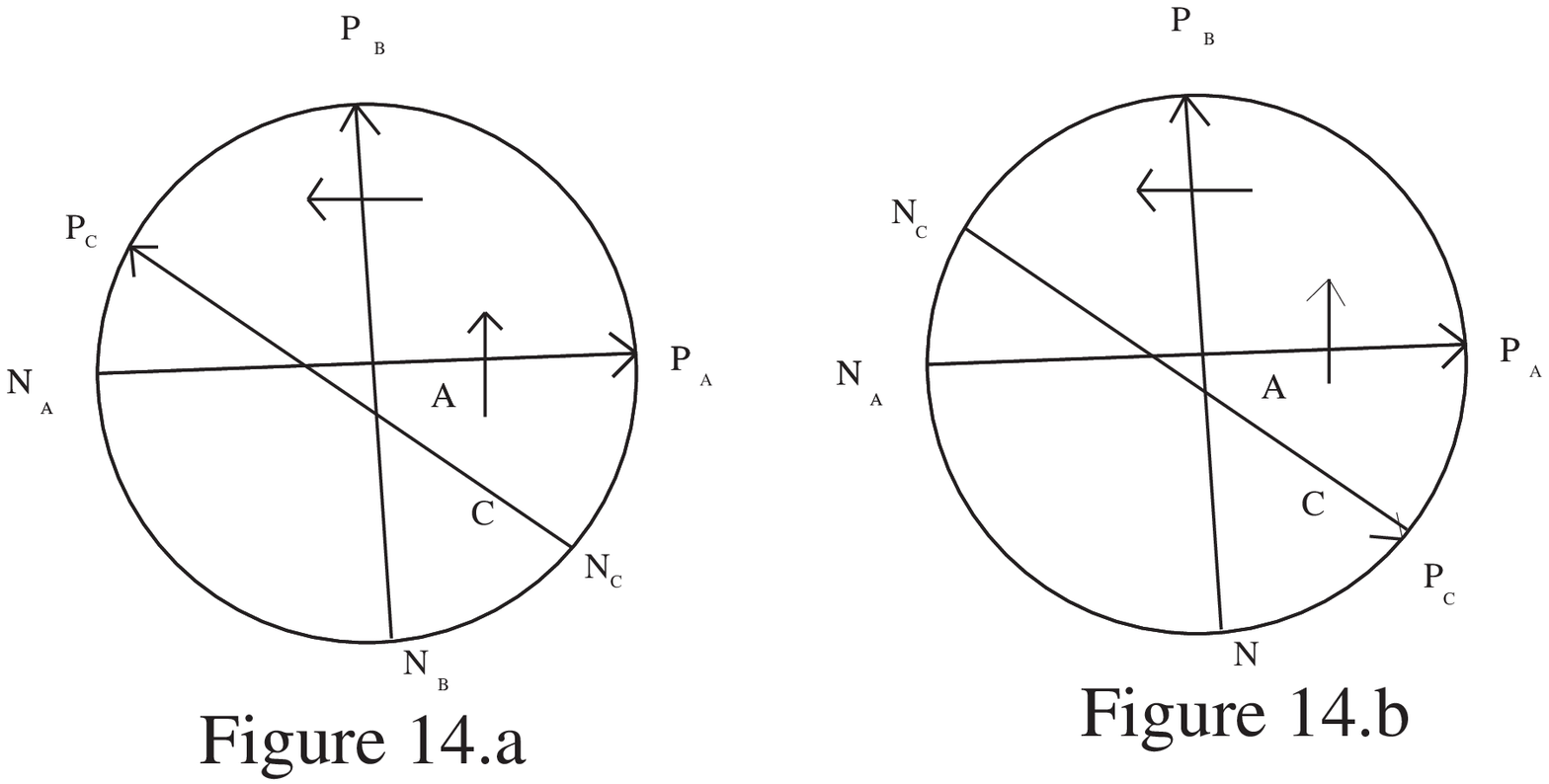 scaled 800}}

We can always change $g$ by an element of $Stab^0B$ to achieve this, by
using Lemma 8.7 above. We call a canonical triple a {\bf good canonical
triple} if in addition $C$ is disjoint from $A$ (See Figures 15a and 15b).

\centerline {\BoxedEPSF{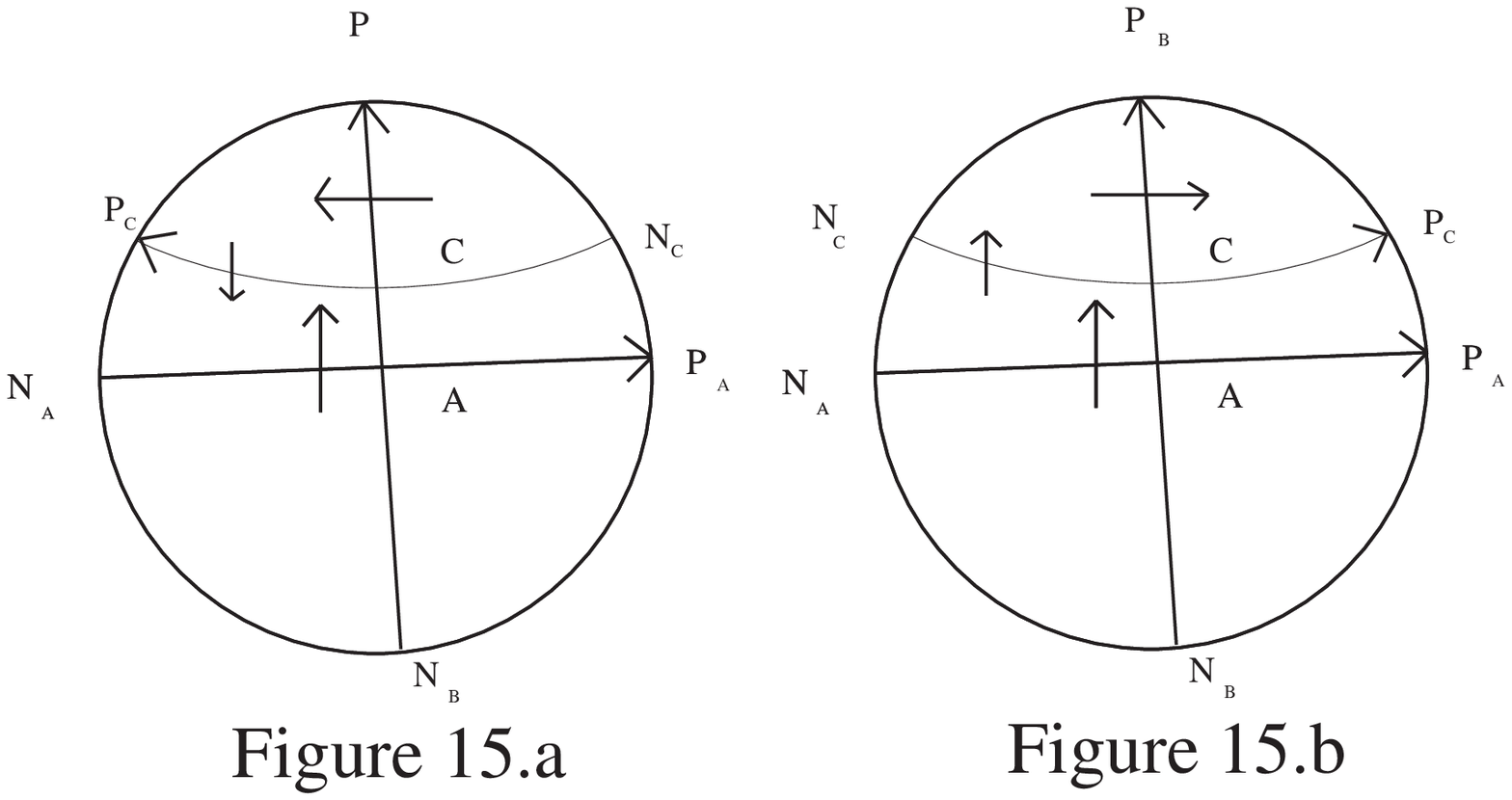 scaled 800}}

We will show that if an axis crosses its translates, then we can obtain good
canonical triples starting with that axis. To handle the o.r. case, we need
a lemma which will be proved in the next section.

\begin{lem}
If $gA$ crosses $A$ for some $g$, then there is a good canonical triple $%
(A,D,E)$.
\end{lem}

{\noindent {\it Proof: }} We proceed much as in Tukia~\cite{tukia:conj}. Let 
$B=gA$ such that the triple $(A,B,gB)$ is a canonical triple. If $gA$
crosses $A$ in an o.r. fashion, then Lemma 9.3 in the next section implies
that $g$ is o.r. and that an axis for $g$ must cross $A.$ Now Corollary 8.4
implies that $g^2A$ and $A$ must be disjoint, so that $(A,B,gB)$ is good,
and the lemma follows.

Otherwise $gA$ crosses $A$ in an o.p. fashion and we consider translates $D$
of $A$ such that $N_D$ lies in $[N_B,P_A]$ and $P_D$ lies in $[N_A,P_B].$
Note that $gB$ is one such. Lemma 8.7 tells us that we can take such a $D$
for which $N_D$ is nearest to $P_A$. This means that there is $h$ in $G$
such that $D=hA$ and $(A,D=hA,E=h^2A)$ is a canonical triple. If $h$ is
o.p.,then $P_E\in [N_A,P_D]\subset [N_A,P_B]$ and by the choice of $D$, $N_E$
cannot be in $[N_D,P_A]$. Thus $(A,D,E)$ is a good canonical triple. If $h$
is o.r. then as in the preceding paragraph, $E$ cannot cross $A$. Thus we
again have a good canonical triple.

In the above situation if $g_A$ is o.r., then if necessary we can modify $h$
by $g_A$ to obtain a new $h$ which is o.r. Thus:

\begin{cor}
With the notation of the above lemma, if $g_A$ is o.r. we can obtain a good
canonical triple $(A,D=hA,hD)$ for which $h$ is o.r.
\end{cor}

\section{Construction of axes using good canonical triples}

We now describe a combinatorial version of a crucial construction of Tukia~%
\cite{tukia:conj}. Tukia showed how given an element $g$ of $G$ such that $%
e(X,<g>)=2$ but $G$ does not split over $<g>$ one can find a ``better''
element. His improvement is in terms of intersection numbers. We could have
adopted his procedure; instead we use shortest tracks and the least area
ideas of Freedman, Hass, and Scott,~\cite{fhs:least} and \cite{fhs:shortest}%
, to obtain a splitting with the first choice of $g$. The proof that this
works is based on Tukia's idea of improving intersection numbers. The
construction is also needed to show that if $hA$ crosses $A$ in an o.r.
fashion, then $h$ is o.r. To have clear pictures we want axes to be minimal
(i.e. covering a shortest finite track) and by the results of \S 7 , we can
do this in the case when the axis stabiliser is o.p. Recall from Lemma 7.10,
that minimal axes in $X$ intersect like lines. This means that if $A,B$ are
minimal axes, then they are disjoint, or coincide, or cross. Further if they
cross then they intersect in a finite graph and both $A-(A\cap B)$ and $%
B-(A\cap B)$ consist of exactly two infinite components each.

We now start with an essential axis $A$ and let $g_A$ be a generator of $%
Stab^0A$. We assume that $A$ is minimal. Consider the situation when $%
(A,gA,g^2A)$ =$(A,B,C)$ is similar to a good canonical triple in the sense
that $A,gA$ cross but $A,g^2A$ do not, but we do not assume that this triple
is actually canonical. If $g$ is o.p. with transverse orientations as shown
in Fig. 12a, then $P_C\in [N_A,P_B]=I_2$ and $N_C\in [P_A,N_B]=I_1$. Recall
from Lemma 8.5 that the fixed points of $g$ lie in $I_2\cup I_4.$ However,
no fixed point can lie in $I_4$ as all points of $I_4$ lie on the negative
side of $A$ and the positive side of $C.$ Thus both the fixed points of $g$
must lie in $I_2$ and in particular they lie on the positive side of $A$ and
all its translates. Denote $g^iA$ by $A_i$; so that $A=A_0,B=A_1,C=A_2$. We
will show now that $A$ intersects $A_i$ if and only if $i$ equals $0,1$ or $%
-1.$ For otherwise, let $k$ be the least integer greater than $1$ such that $%
A$ intersects $A_k.$ Then the set of all points in $X$ which lie on the
positive side of each of $A_0,\ldots ,A_k$ is compact, contradicting the
fact that both the fixed points of $g$ lie on the positive side of $A$ and
all its translates. This allows us to construct an axis for $g$ as follows.
By the lemmas above, one of the two infinite components of $A-A\cap B$
together with one of the infinite components of $B-A\cap B$ bounds the
infinite component of $X-(A\cup B)$ which contains $I_2$ in its boundary. We
denote this pattern by $A\#B=A_0\#A_1$. Inductively we construct (See Figure
16) $A_{-n}\#...\#A_0\#...\#A_n=L_n$, say. Since the collection $g^nA$ is
locally finite, we obtain a limit pattern $A_g$ of $L_n$ such that $gA_g=A_g$
and the quotient of $A_g$ by the action of $g$ is a finite pattern $t_g$ in $%
X_g$. We claim that some component of $t_g$ is an essential track, so that,
in particular, $<g>$ has two co-ends.

\centerline {\BoxedEPSF{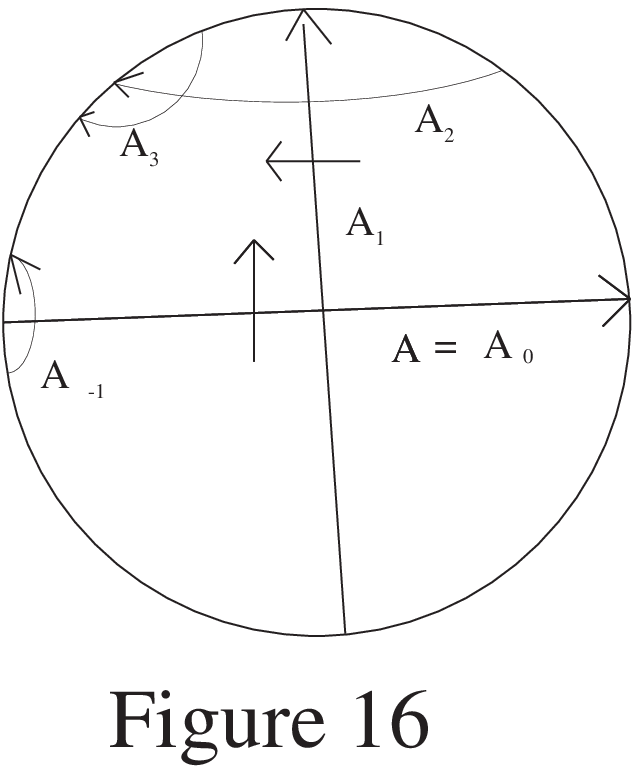 scaled 800}}

Recall that we are assuming that $g$ is o.p. Thus if no component of $t_g$
is essential in $X_g$, then every component of $t_g$ bounds a compact subset
of $X_g$. If we remove these compact sets from $X_g$ and take the closure $%
Y_g,$ then $Y_g$ carries $<g>.$ This is because if a component of $t_g$
carries the trivial group, then it bounds a compact subset of $X_g$ which
also carries the trivial group, by Lemma 3.1. Thus the pre-image of $Y_g$ in 
$X$ is its universal cover $Y,$ and the number of ends of $Y$ and of $X$ are
the same. Recall that the pattern $A_g$ intersects $A$ in the part of $A$
between $A_1$ and $A_{-1}.$ Now take a large enough $N$ so that $g_A^NA_g$
is to the right of $B$. Then the portion of $A$ between $B$ and $g^NB$ is
compact and divides $Y$ into two infinite parts, which contradicts the
assumption that $X$ and hence $Y$ has only one end. Thus we have:

\begin{lem}
If $A$ is a minimal essential axis, if $g$ is o.p. and if $A,gA$ intersect
but $A,g^2A$ do not, then $<g>$ has two co-ends.
\end{lem}

This is one place where the use of minimal axes seems convenient even though
the lemma itself can be formulated purely algebraically. It should be
possible to prove this lemma by more algebraic means and avoid minimal axes
completely and proceed with intersection numbers as in Tukia~\cite
{tukia:conj}. The case of o.r. elements is somewhat easier.

\begin{lem}
If $A$ is a minimal essential axis, if $A$ and $gA$ have o.r. crossing, and
if $A,g^2A$ do not cross, then $g$ is o.r. and has fixed points in $I_1$ and 
$I_3,$ one in each.
\end{lem}

The construction is similar, but in this case the two sides of the boundary
into which $A_g$ separates $\partial X$ are clearly non-empty, as they are $%
\bigcup g^{2i}[N_{-1},P_A]$ and $\bigcup g^{2i+1}[N_{-1},P_A]$ respectively,
where $P_{-1}$ and $N_{-1}$denote the positive and negative endpoints of $%
A_{-1}.$ See Figure 17.

\centerline {\BoxedEPSF{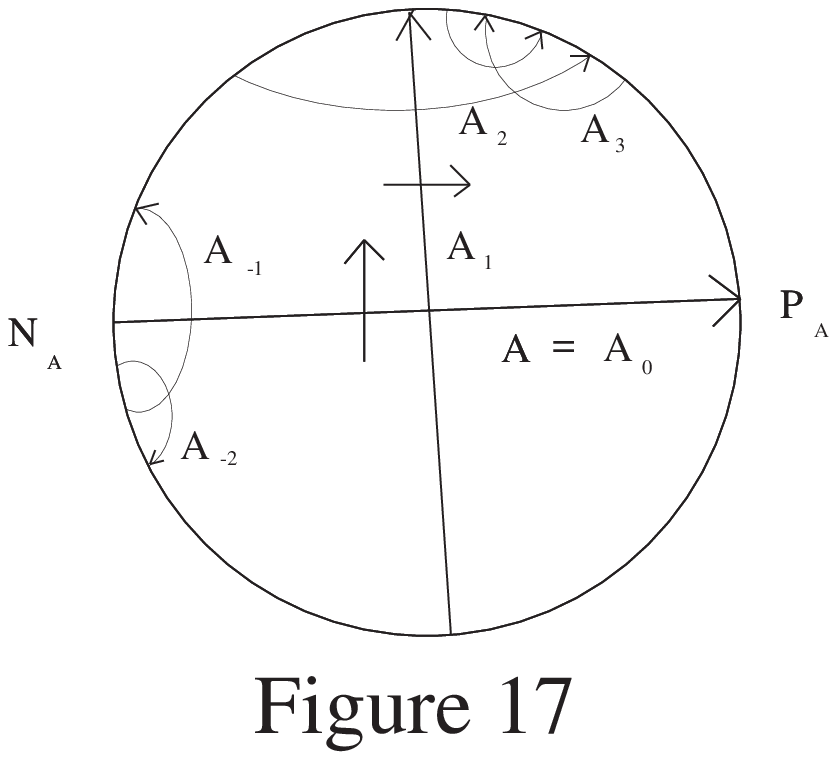 scaled 800}}

We use this to prove;

\begin{lem}
If $hA$ crosses $A$ in an o.r. fashion, then $h$ is o.r.
\end{lem}

{\noindent {\it Proof: }} We will assume that $A$ is minimal. As the fixed
points of $h$ in $\partial X$ must be distinct from the end points of $A,$
there is a greatest power $h^n$ of $h$ such that $h^nA$ crosses $A.$ If $n$
is even, then this crossing must be o.p. As $h^{2n}A$ does not cross $A,$
Lemma 9.1 shows that $h^n$ has an essential axis. It follows that the same
is true for $h.$ Now it follows that $h$ must be o.r. If $n$ is odd, then
the crossing of $A$ and $h^nA$ must be o.r. and Lemma 9.2 shows that $h^n$
is o.r. which also implies that $h$ is o.r. This completes the proof of
Lemma 9.3.

We will need the following consequence in the next section. This easily
follows from the above since $h$ and hence $h^2$ has an essential axis.

\begin{cor}
If $A$ is a minimal essential axis and $hA$ crosses $A$ in an o.r. fashion,
then $h^2A$ is disjoint from $A$.
\end{cor}

Whether $A$ is minimal or not, as long as it is essential, by replacing $A$
by a minimal axis with the same end points we conclude:

\begin{cor}
If $A$ is an essential axis and $hA$ crosses $A$ in an o.r. fashion, then
the end points of $h^2A$ are on the same side of $A$.
\end{cor}

\section{Proof of the annulus theorem in the surface type case}

We start on the proof of the annulus theorem in the surface type case. We
first consider the case when $G$ has o.r. elements and prove the following
special case of the annulus theorem.

\begin{lem}
Let $G$ be a torsion free hyperbolic group with one end and suppose that $G$
is of surface type. If $G$ contains an o.r. element, then there is an o.r.
element $g$ of $G$ such that $G$ splits over $<g^2>$ as an amalgamated free
product with $<g>$ as one of the two factor groups.
\end{lem}

{\noindent {\it Proof: }}We choose $g$ to be an o.r. element of $G$ such
that $c(g)$ is minimal, where $c(g)$ is the complexity of a shortest
one-sided track $s$ in $X_g.$ Note that $g$ will be indivisible. Let $S$
denote the corresponding axis of $g.$ Also let $t$ be a shortest essential
track in $X_{g^2},$ and recall from Lemma 7.6 that $2c(g)\geq c(g^2)$ where $%
c(g^2)$ is the complexity of $t$. Let $A$ denote the axis of $g^2$
corresponding to $t.$ Note that $S$ and $A$ are both stabilised by $g^2$ and
so have the same end points in $\partial X.$ Also $g$ stabilises $S$ but may
not stabilise $A,$ and $A$ is minimal but $S$ need not be minimal. If $S$
and its translates never cross, then Theorem 5.1 shows that $G$ must split
over a subgroup commensurable with $<g>,$ and one can then show that $G$ has
the required splitting over $<g^2>.$ However, as $G$ is of surface type, we
can use the properties of minimal axes to give a much simpler proof. As $S$
and its translates never cross, the same holds for $A$ and its translates.
As $A$ is minimal, Lemma 7.10 shows that $A$ is disjoint from or coincides
with each translate. It follows that $A$ covers a finite track $u$ in $X_G$
which must carry $<g>$ or $<g^2>.$ If $u$ carries $<g>$ then $g$ stabilises $%
A.$ In this case, the projection of $A$ into $X_g$ must be a one-sided track
so that $u$ also must be one-sided. We replace $u$ by the boundary of a
regular neighbourhood in $X_G.$ This will be a two-sided track $v$ carrying $%
<g^2>$ and $v$ separates $X_G$ into two pieces one of which is a regular
neighbourhood of $u$ and so carries $<g>.$ It follows at once that $G$ has
the required splitting. Now suppose that $u$ carries $<g^2>.$ Thus $u$ lifts
to $X_{g^2}$ and projects to an embedded two-sided track in $X_g$ which we
denote by $u^{\prime }.$ The region of $X_{g^2}$ between $u$ and its
translate under the covering involution projects to a compact connected
subset $Z$ of $X_g$ which carries $<g>$ and is bounded by $u^{\prime }.$ We
claim that $Z$ projects into $X_G$ by a homeomorphism so that $u$ bounds a
copy of $Z$ in $X_G.$ Again this implies that $G$ has the required
splitting. To prove the claim, let $p:X_g\rightarrow X_G$ denote the
covering projection and consider $p^{-1}(u)\cap Z.$ This consists of $%
u^{\prime }$ and possibly other finite covers of $u,$ contained completely
in $Z.$ Let $v$ denote a finite cover of $u$ contained in $Z.$ It must carry 
$<g^n>$ for some positive $n$ and its pre-image in $X$ must be a translate $%
kA$ of $A$ which is stabilised by $<g^n>.$ It follows that $kg^2k^{-1}$
equals $g^n$ or $g^{-n}.$ This is only possible if $n=2$ and $k$ is a power
of $g,$ as $G$ is word hyperbolic. This implies that $v$ equals $u^{\prime
}, $ so that $p^{-1}(u)\cap Z=u^{\prime }.$ Let $\Sigma $ denote the subset $%
\{z\in Z:\exists y\in Z,y\neq z,p(y)=p(z)\}$ of $Z.$ Then $\Sigma $ must be
disjoint from $u^{\prime }.$ It follows that $\Sigma $ is open and closed in 
$Z,$ and hence equals $Z$ or is empty. The fact that $\Sigma $ is disjoint
from $u^{\prime }$ implies that $\Sigma $ must be empty, so that $Z$ must
project into $X_G$ by a homeomorphism as claimed.

Now we consider the case where some translate of $A$ crosses $A.$ We will
also assume that the intersection of any two translates of $A$ is
transverse. The general case can again be handled by using the Meeks-Yau
trick. Some translate of $S$ must also cross $S,$ so that Corollary 8.10
tells us that there is a good canonical triple $(S,hS,h^2S)$ with $h$ o.r.
It follows that there is a good canonical triple $(A,D=hA,hD)$ with $h$ o.r.
This means that there is no translate $E$ of $A$ having $P_E\in (N_D,P_A)$
and $N_E\in (N_A,P_D)$. In particular the axes $(g_D^{2m})hD$ do not cross $%
A,$ for any value of $m,$ where $g_D=hgh^{-1}.$ Choose $m$ to be the least
value such that the end points of $(g_D^{2m})hD$ lie above $A$ and replace $%
h $ by $(g_D^{2m})h.$ This does not change $D$, so we obtain a new triple $%
(A,D=hA,hD)$ and $hD$ still does not cross $A.$ We have the end points of $%
hD $ above $A$ and now the end points of $g_D^{-2}hD$ lie below $A$ (see
Figure 18).

\centerline {\BoxedEPSF{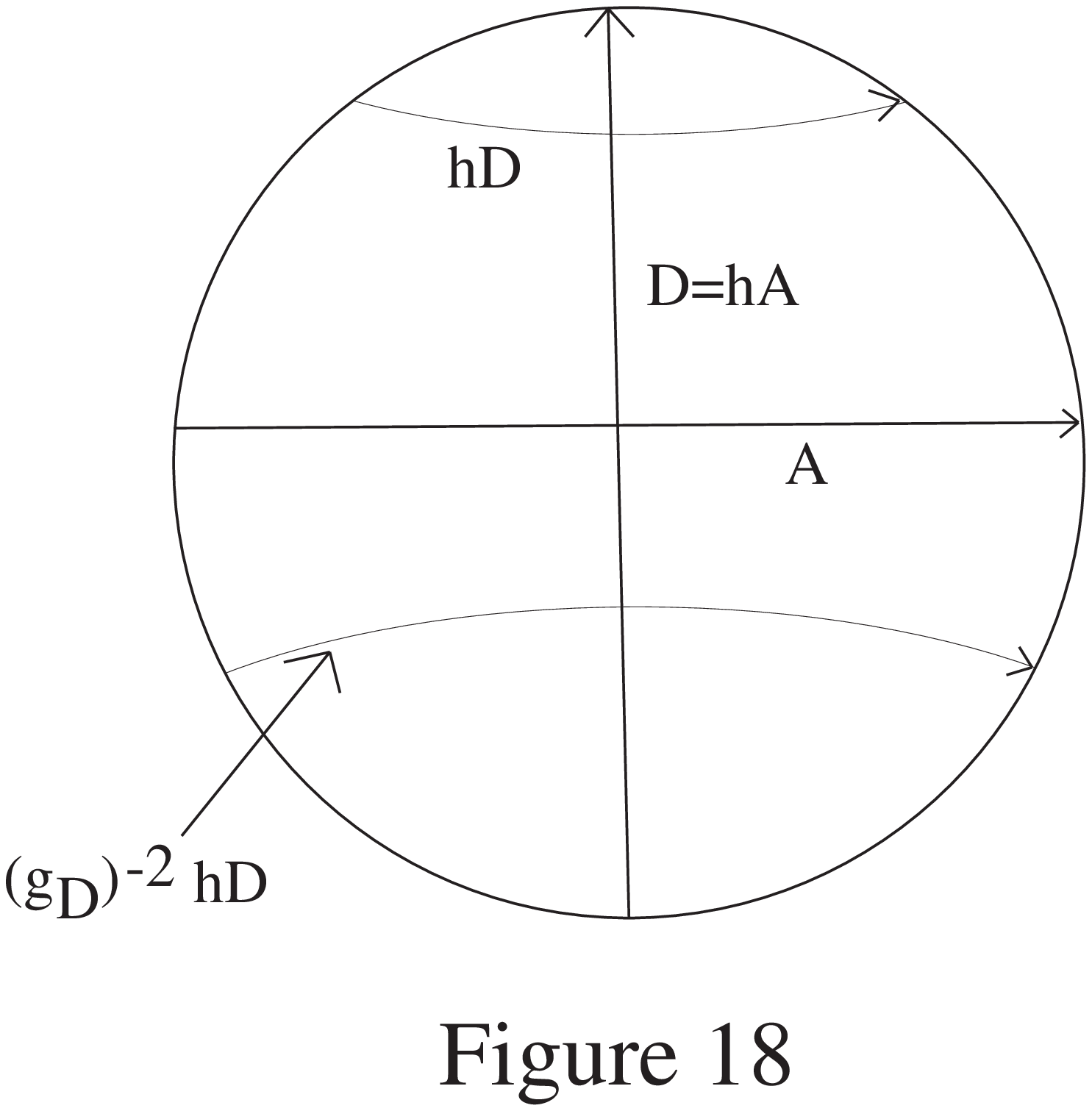 scaled 400}}

Let $k$ denote $g_D^{-2}h$ and consider the triples $(A,D,hD)$ and $(kD,D,A)$%
. The second triple is also like a good canonical triple since $k^{-1}D=A,$
but it need not be canonical. Also $kD$ and $D$ cross o.r. so that $k$ must
be o.r. If $D_1$ denotes the portion of $D$ between $A$ and $hD$ and $D_2$
the portion of $D$ between $g_D^{-2}hD$ and $A$, we have $%
c(D_1)+c(D_2)=c(g_D^2)=c(g^2)$. Note that the proof of Lemma 9.2 implies
that the image of $D_1$ in $X_h$ contains a one-sided track carrying $<h>,$
and the image of $D_2$ in $X_k$ contains a one-sided track carrying $<k>.$
Thus either $c(h)$ or $c(k)$ is, by rounding off corners, strictly less than
half of $c(g^2)$. This contradicts our choice of $c(g)$. It follows that the
translates of $A$ cannot cross, so that $G$ splits over $<g^2>$ as claimed.
This completes the proof of the Annulus Theorem when $G$ contains o.r.
elements.

In order to complete our proof of the Annulus Theorem, we consider the case
when $G$ has no orientation reversing elements. In this case we choose an
element $g$ of $G$ such that $c(g)$ is minimal among all elements of $G$
with two co-ends, and let $A$ denote the corresponding minimal essential
axis for $g.$ If $hA$ intersects $A$ for some $h$, then by the results of \S
8, we can choose $h$ so that $(A,hA,h^2A)$ is a good canonical triple (See
Figure 19. This is similar to Figure 5 in Tukia~\cite{tukia:conj}).

\centerline {\BoxedEPSF{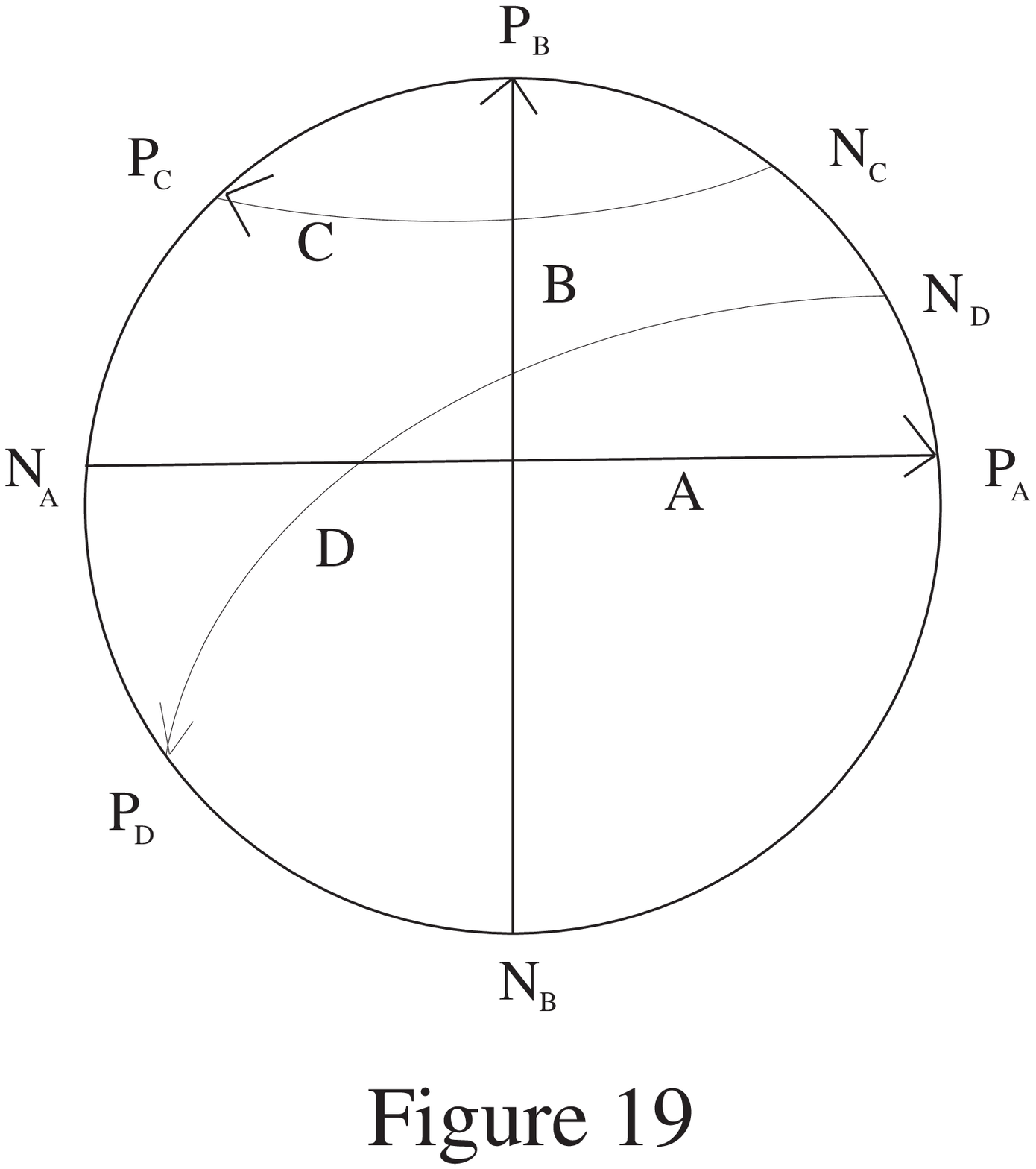 scaled 400}}

Note that $h$ is automatically o.p. Now consider $g_B^{-1}C=D$. Since $%
(A,B,C)$ is a good canonical triple $g_B^{-1}P_C=P_D$ is automatically below 
$A$. If $N_D$ is also below $A,$ we can argue as in the proof of Lemma 10.1,
that either $c(h)$ or $c(g_B^{-1}h)$ is strictly less than $c(g),$
contradicting the choice of $g.$ If $N_D$ is above $A$ we arrive at a
contradiction as follows. Consider the fixed points of $g_B^{-1}h$. Since $%
B=hA=g_B^{-1}hA$, the fixed points of $g_B^{-1}h$ should be in $%
[N_A,P_B]\cup [N_B,P_A],$ by Lemma 8.5. Similarly, since $%
D=g_B^{-1}C=g_B^{-1}hB$, the fixed points of $g_B^{-1}h$ should be in $%
[P_D,N_B]\cup [N_D,P_B]$. But these sets are disjoint and this contradiction
completes the proof of the annulus theorem in the surface type case.

\clearpage
\obeylines
\parindent0pt
Peter Scott
Mathematics Department, University of Michigan, Ann Arbor, MI 48109
\quad pscott@umich.edu

\bigbreak

Gadde A. Swarup
Mathematics Department, University of Melbourne, Parkville, Victoria 3052,
\quad  Australia
\quad gadde@maths.mu.oz.au 


\begin{thebibliography}{99}
\bibitem{bf:comb}  M. Bestvina and M. Feighn. {\em A combination theorem for
negatively curved groups} Journal of Differential Geometry,35, 85-101,(1992).

\bibitem{bo:cut}  B. Bowditch. {\em Cut points and canonical splittings of
hyperbolic groups} Preprint, University of Southampton, 1995.

\bibitem{bo:tree}  B. Bowditch. {\em Group actions on trees and dendrons}
Preprint, University of Southampton, 1995.

\bibitem{coo:thesis}  M. Coornaert {\em Sur les groupes properement
discontinus d'isometries des espaces hyperboliques au sens de Gromov}
Thesis, Universite Louis Pasteur, Strasbourg (1990).

\bibitem{mjd:acc1}  M. J. Dunwoody. {\em Accessibility and groups of
cohomological dimension one} Proc. London Math. Soc. 38, 193-215, (1979)

\bibitem{mjd:acc2}  M. J. Dunwoody. {\em The accessibility of finitely
presented groups} Inv. Math. 81, 449-457(1985).

\bibitem{ef:conv}  E. Freden. {\em Negatively curved groups have convergence
property} Ann.Acad.Sci.Fenn.Ser.A Math. 20, 333-348 (1995)

\bibitem{fhs:least}  M. H. Freedman, J. Hass and G. P. Scott. {\em Least
area incompressible surfaces in 3-manifolds} Inv. Math 71, 609-642 (1983)

\bibitem{fhs:shortest}  M. H. Freedman, J. Hass and G. P. Scott. Closed
Geodesics on Surfaces, Bull. London Math. Soc. 14(1982), 385-391.

\bibitem{gabai:co}  D. Gabai {\em Convergence groups are Fuchsian groups}
Annals of Maths. 136, 447-510 (1992)

\bibitem{gh:book}  E. Ghys and P. de la Harpe. {\em Sur les groupes
hyperboliques d'apres Mikhail Gromov} Progress in Maths. No. 83, Birkhauser
(1990).

\bibitem{jr:pl}  W. Jaco and J. H. Rubinstein. {\em PL minimal surfaces in
3-manifolds} J. Differential Geometry 27, 493-524 (1988).

\bibitem{pau:outer}  F. Paulin {\em Outer automorphisms of hyperbolic groups
and small actions on R-trees} in Arboreal Group Theory(edited by R. Alperin)
MSRI Publications 18, 331-343, (1991)

\bibitem{rs:rigid}  E. Rips and Z. Sela. {\em Structure and rigidity of
hyperbolic groups} Geom. and Funct. Analysis 4, 337-371 (1994)

\bibitem{sc1:newpro}  G. P. Scott. {\em A new proof of the annulus and torus
theorems} American J. of Math. 102, 241-277 (1980).

\bibitem{sc2:lerf}  G. P. Scott. {\em Subgroups of surface groups are almost
geometric} Journal of L.M.S., 17 , 555-565 (1978)

\bibitem{sc3:ends}  G. P. Scott. {\em Ends of pairs of groups} J. of Pure
and Applied Algebra,11,179-198 (1977)

\bibitem{sela:rigid2}  Z. Sela. {\em Structure and rigidity in (Gromov)
hyperbolic groups} Annals of Math. (1995)

\bibitem{st:spl}  J. R. Stallings. {\em Group theory and 3-dimensional
manifolds} Yale University Press , New Haven (1971)

\bibitem{sw:ends2}  G.A. Swarup {\em On ends of pairs of groups} J. of Pure
and Applied Algebra,11,179-198 (1993)

\bibitem{sw:access}  G.A. Swarup {\em A note on accessibility} in Geometric
Group Theory, ( Edited by G.Niblo and M.Roller) Cambridge University Press,
204-207, (1993)

\bibitem{sw:cutpoint}  G.A. Swarup {\em On the cut point conjecture},
Preprint (1996)

\bibitem{tukia:conj}  P. Tukia. {\em Homeomorphic conjugates of Fuchsian
groups} J. reine angew. Math. 391, 1-54 (1988).

\bibitem{tukia:conv}  P. Tukia {\em Convergence groups and Gromov's metric
hyperbolic spaces} New Zealand J. Math 23, 157-187 (1994)
\end{thebibliography}
\end{document}